# Forward-looking evolutionary game dynamics subject to exploration cost

Hidekazu Yoshioka[1, *]

[1] Graduate School of Advanced Science and Technology, Japan Advanced Institute of Science and Technology, 1-1 Asahidai, Nomi, Ishikawa, Japan

[*] Corresponding author: yoshih@jaist.ac.jp, ORCID: 0000-0002-5293-3246

***Abstract***

We extend classical evolutionary game dynamics based on the momentary action choices of agents by accounting for two elements: forward-looking behavior and exploration cost. We focus on pairwise comparison protocols that cover major evolutionary game dynamics, such as replicator and logit models. In the proposed mathematical framework, agents update their actions by paying a cost so that a utility or its relative difference is maximized. We show that forward-looking behavior can be modeled as a coupling between the evolutionary game dynamic and static Hamilton–Jacobi–Bellman equation: a mean field game. The exploration cost and its constraint are naturally related to these equations as a function of the optimal Lagrangian multiplier serving as a relaxation parameter, and it is incorporated into the game as a constraint. We show that under certain conditions, our evolutionary game dynamic admits a unique solution. Finally, we computationally investigate one- and two-dimensional problems.

***Statements & Declarations***

**Competing interests:** The author has no relevant financial or nonfinancial interests.

**Data availability:** Data will be made available upon reasonable request to the author.

**Funding:** Japan Society for the Promotion of Science (KAKENHI, No. 25K00240, 26K00357)

**Acknowledgments:** NA

**Declaration of Generative AI and AI-assisted technologies in the writing process:** The author used generative AI only for proofreading the manuscript.

# 1. Introduction

## 1.1 Study background

Evolutionary game dynamics (EGDs) describe the actions of interacting agents aiming at maximizing their utilities [21] and have been used as mathematical models for social phenomena. The governing equation of an EGD is typically given by a difference or differential equation, and its structure is problem dependent. Recent examples of EGDs in applications include the difference equation model for the management of public goods [26], difference equation model for the evolution of cooperation [42], ordinary differential equation (ODE) model for social norms in social-ecological systems [40], ODE model for epidemic control [15], and partial integro-differential equation (PIDE) model for social dilemmas [47].

In this paper, we focus on continuous-time and continuous-state EGDs, which are PIDEs, because of their versatility for describing social phenomena, and we can obtain difference equation models or ODE models by suitable discretization. EGDs as PIDEs have been theoretically studied in Cheung [12] in the following form (more rigorously presented later):

$$\underbrace{\frac{\mathrm{d}\mu_t(\mathrm{d}x)}{\mathrm{d}t}}_{\text{Change of the probability of agent actions}} = \underbrace{\int_{z\in\Omega}\rho\big(U(z,\mu_t),U(x,\mu_t)\big)\lambda_t(\mathrm{d}x)\mu_t(\mathrm{d}z)}_{\text{Inflow of probability}} - \underbrace{\int_{z\in\Omega}\rho\big(U(x,\mu_t),U(z,\mu_t)\big)\lambda_t(\mathrm{d}z)\mu_t(\mathrm{d}x)}_{\text{Outflow of probability}}, \quad (1)$$

where $\mu_t$ is the probability measure of agent actions $x$ in an action space $\Omega$ at time $t>0$, $\rho$ is a bivariate function called the protocol in the sequel, $U(x,\mu_t)$ is the utility of agents who choose the action $x$, and $\lambda_t$ is some probability measure that may depend on $\mu_t$. The EGD (1) is nonlinear because the utility $U$ depends on the solution $\mu_t$. A fundamental existence and uniqueness result of this EGD has been established in Theorem 1 in Cheung [12]. An important observation of (1) is that the evolution of agent actions, given a utility $U$ and a probability measure $\lambda$, is determined by the protocol $\rho$, which typically depends on the utility difference $U(x,\mu_t)-U(z,\mu_t)$ (difference maximization case) or the utility $U$ itself (value maximization case). Well-studied protocols for the difference maximization case include the replicator [38] and Brown–von Neumann–Nash (BNN) models [20], and those for the value maximization case include the logit model [31]. These protocols have been investigated not only theoretically but also from application-oriented standpoints such as decision-making with waiting times [18,43], dissipativity and stationary states focusing on specific utilities [2], auctions [34], pollution and environmental regulations [30,41], and strategic decision-making to attract and retain talented employees by organizations [14].

At least two underexplored research challenges concerning EGDs exist: forward-looking behavior and exploration cost. First, the utility $U$ typically depends on the current state of agent actions but not their future predictions, even though it is natural for agents to make current decisions while predicting future events (called forward-looking behavior), as in typical dynamic optimization problems such as dynamic programming [7]. In this context, mean-field game (MFG) models are natural extensions of (1), where an EGD as an initial-value problem for updating actions is coupled with a Hamilton–Jacobi–Bellman (HJB) equation as a terminal-value problem for future predictions, forming a forward-backward system. A practical drawback of this approach, which is common to other MFGs, is that the well-posedness and computability of the resulting system are often difficult or unresolved because of the forward-backward nature [49,50]. Couplings between a static EGD and a static HJB equation also exist, but they lose the dynamic nature of EGDs [6,51]. Recently, Iijima and Oyama [23] proposed coupling an EGD with a "static" and hence time-independent uncontrolled HJB equation, which is free from the forward-backward nature; however, how the resulting system behaves has not been thoroughly investigated.

The second issue is the exploration cost, meaning that the EGD (1) does not involve a mechanism such that making a better (e.g., more accurate and/or fast) decision would be more costly. This point has been modeled through information-processing costs in the context of optimal control problems [1,48,53]. In this framework, it is assumed that there exists a maximum cost to be paid for processing the information necessary for making decisions and that relaxing this constraint probably yields better results for the decision maker. By applying this idea to the EGD (1), the protocol $\rho$ can be determined through a constrained optimization problem to maximize some outcome (e.g., increasing the utility more efficiently, faster convergence to an equilibrium, etc.) subject to the cost of information processing. Moreover, once an evaluation method of the cost is established, it will become possible to compare the performance among different protocols for each prescribed cost value in the future. To the best of the author's knowledge, such approaches have not been examined thus far.

In summary, conventional studies related to EGD (1) do not thoroughly investigate or account for forward-looking behavior or exploration cost despite both being important elements in the advanced modeling of social phenomena; thus, this research gap remains to be resolved.

### 1.2 Aim and contribution

Motivated by the background reviewed above, this study aims at the formulation, computation, and application of an EGD accounting for both forward-looking behavior and exploration cost. Our contributions to achieving this aim are explained as follows. Note that our contribution is novel even if there is no forward-looking behavior or no exploration cosin EGDs.

First, we formulate an extended EGD (called forward-looking EGD in this paper) under a general condition where the utility $U$ is replaced by a forward-looking version, called the value function in accordance with MFGs, which is governed by a static HJB equation with the classical utility as a source term and a discount term to modulate the myopicity of agents. This method is qualitatively analogous to that in Iijima and Oyama [23] if the HJB equation is linearized. Another key element, which is the exploration cost, is then formulated by using a constrained optimization problem whose solution is the protocol $\rho$. The exploration cost is incorporated into the HJB equation with a Lagrangian-multiplier approach. The optimal Lagrangian multiplier that satisfies the cost constraint enters the forward-looking

EGD as a time-dependent relaxation parameter, which effectively controls the speed and/or accuracy of evolution of the probability measure $\mu$ of agent actions.

Second, we mathematically study specific cases where the protocol $\rho$ is of the BNN, replicator, or logit type. We show that the functional form of cost for each protocol is found explicitly, while the HJB equation does not admit a closed-form solution except for the logit model. We show that forward-looking EGD yields a unique solution globally in time under suitable regularity conditions of utilities, and a regularization term for BNN and replicator models. This part is based on a Banach fixed-point approach for dealing with the HJB part in conjunction with the unique existence criterion for EGDs (Theorem 1 in Cheung [12]). An additional difficulty compared with classical EGDs is that our dynamic is a coupling between the PIDE of the form (1) and the static HJB equation with a cost constraint and hence is more strongly nonlinear than the classical ones are. We present an explicit sufficient condition for the unique existence of solutions to our forward-looking EGD.

Mathematically, forward-looking EGD can be seen as an MFG subject to a cost constraint, which by itself is an unexplored topic in MFGs; indeed, in most cases, an MFG accounts for costs in the objective function by using a weighting coefficient and hence does not explicitly constrain the cost [e.g., 11,16,45]. This point is a difference between existing and our MFG formulations. In past studies, the characterization of an EGD from a nonclassical viewpoint has sometimes been addressed, e.g., Jaćimović [24,25] characterized a replicator dynamic as a gradient flow of optimization problems. Additionally, replicator dynamics were characterized as EGDs under which exponential measures are invariant [33]. BNN, replicator, and logit protocols without forward-looking behavior are characterized by a large-discount limit of the HJB equation under which the value function coincides with utility.

Third, we present a simple numerical method for the discretization of forward-looking EGDs. This numerical method uses a naïve finite difference discretization in conjunction with an iteration scheme for computationally realizing the cost constraint. The discretization scheme itself is based on standard numerical techniques, but importantly, specialized algorithms are not always necessary for computing our EGDs. We use this numerical method to investigate our EGDs with different protocols. We computationally investigate forward-looking EGDs with different protocols based on the numerical method to investigate their performance.

Consequently, this study addresses the formulation, analysis, computation, and application of a new EGD model in a consistent manner.

### 1.3 Structure of this paper

**Section 2** presents classical and forward-looking EGDs. **Section 3** theoretically studies the uniqueness and existence of forward-looking EGDs with specific protocols. **Section 3** also presents a numerical algorithm to discretize forward-looking EGDs. **Section 4** computationally demonstrates forward-looking EGDs. **Section 5** summarizes this study and its future perspectives. **Appendix** presents auxiliary results; **Section A.1** defines Nash equilibria, **Section A.2** presents proofs of the propositions in the main text, **Section A.3** presents a two-dimensional numerical method, and **Section A.4** presents a convergence study of the numerical method. **Section A.5** explains an analytical result about the computational results in **Section 4**.

## 2. Evolutionary game dynamics

Throughout this paper, we consider EGDs with a one-dimensional compact action space, which is a typical setting for studying EGDs with continuous actions [e.g., 20,31,38]. The domain of agent actions is set as $\Omega = [0,1]$ without any loss of generality. Time is a continuous parameter and is denoted by $t \geq 0$. The space of all probability measures on $\Omega$ is denoted by $\mathcal{P}(\Omega)$. The probability measure of agent actions at time $t$ is denoted by $\mu_t$. The collection of all Borel measurable sets in $\Omega$ is denoted by $\mathcal{B}(\Omega)$. The one-dimensional real space is denoted by $\mathbb{R}$, and its nonnegative part by $\mathbb{R}_+$. For any $a,b \in \mathbb{R}$, we write $\max\{a,b\}$ as $a \vee b$ and $\min\{a,b\}$ as $a \wedge b$. For any $a \in \mathbb{R}$, we write $\max\{a,0\} = (a)_+$.

### 2.1 Protocols

We present a novel view of the protocol $\rho$. The EGD, which was formally presented in (1), is more rigorously presented as follows because $\mu$ is not necessarily a function:

$$\begin{aligned}\frac{\mathrm{d}\mu_t(A)}{\mathrm{d}t} &= \int_{x\in A}\int_{z\in\Omega}\rho\left(U(z,\mu_t),U(x,\mu_t)\right)\lambda_t(\mathrm{d}x)\mu_t(\mathrm{d}z)\\ &\quad-\int_{x\in A}\int_{z\in\Omega}\rho\left(U(x,\mu_t),U(z,\mu_t)\right)\lambda_t(\mathrm{d}z)\mu_t(\mathrm{d}x)\end{aligned} \tag{2}$$

for any Borel measurable set $A \subset \mathcal{B}(\Omega)$ and time $t > 0$. The EGD (2) is subject to an initial condition $\mu_0 \in \mathcal{P}(\Omega)$. Utility $U:\Omega\times\mathcal{P}(\Omega)\to\mathbb{R}$ to be maximized by agents is a mapping from the action and its distribution to a real space.

Typical protocols $\rho$ depend on the utility difference $U(x,\mu_t)-U(z,\mu_t)$ or the current utility $U$ [12]; hence, with an abuse of notations, we rewrite the bivariate protocol $\rho\left(U(z,\mu_t),U(x,\mu_t)\right)$ as the univariate version $\rho\left(U(x,\mu_t)-U(z,\mu_t)\right)$ for the difference maximization case (e.g., BNN and replicator models) and $\rho\left(U(z,\mu_t)\right)$ for the value maximization case (e.g., logit model). Because agents try to maximize their utilities or their difference, it can be natural to assume that $\rho$ is an increasing and nonnegative function. For the replicator model, we have $\rho(v)=(v)_+$ with $\lambda_t=\mu_t$, and for the BNN model, we have $\rho(v)=(v)_+$ ($v\in\mathbb{R}$) with $\lambda_t$ being a probability measure whose support is $\Omega$. Finally, for the logit model, we set $\rho\left(U(z,\mu_t)\right)=e^{\eta^{-1}U(z,\mu_t)}/\int_{\Omega}e^{\eta^{-1}U(y,\mu_t)}\mathrm{d}y$ with some $\eta>0$. For the logit model, for simplicity, we assume a uniform distribution ($\lambda_t(\mathrm{d}x)=\mathrm{d}x$) in the rest of this paper. In **Section 2.2**, the form of each $\rho$ is characterized through optimization problems.

## 2.2 An optimization view of classical EGDs: cost constraints

### 2.2.1 General setting

We characterize $\rho$ as a maximizer of a constrained optimization problem that represents an exploration (or updating) strategy of agent actions for reinterpretation from the perspective of cost and effectiveness. The inverse-control approach by Yoshioka [49] provides a hint to achieve this objective. Consider the following optimization problem for each fixed $x \in \Omega$ and $\sigma, \theta \in \mathcal{P}(\Omega)$: Find

$$\begin{gathered} \sup_{u_x} \underbrace{\int_{z\in\Omega} u_x(z)\left(U(z,\sigma)-U(x,\sigma)\right)\theta(\mathrm{d}z)}_{\text{Net utility}} \\ \text{s.t.} \\ \underbrace{\int_{x\in\Omega}\int_{z\in\Omega} F\left(u_x(z)\right)\theta(\mathrm{d}z)\theta(\mathrm{d}x)}_{\text{Cost}} \leq \varepsilon \end{gathered} \tag{3}$$

with a constant $\varepsilon > 0$ and a function $F:\mathbb{R}_+ \to \mathbb{R}_+$. This is a maximization problem of the difference between the current utility $U(x,\sigma)$ and the updated utility $U(z,\sigma)$, and the control variable $u_x$ represents the jump rate from the current state $x$ to the updated state $z$. The domain of $u_x$ depends on protocols as discussed below. The first and last lines in (3) represent the net utility gained by updating the action and a constraint of exploration cost for updating agent actions, which is modeled by using a function $F$. With this formulation, we want to regard the optimal control $u_x = u_x^*$ as $\rho$, so that each protocol arises from a maximization principle subject to a cost constraint. We assume that $F$ is strictly convex and continuously differentiable in the sequel, with which we can cover the BNN, replicator, and logit models that are of interest in this paper. **Table 1** summarizes the settings of the coefficients for the BNN, replicator, and logit models.

**Table 1.** Setting of coefficients for replicator, BNN, and logit models.

<table>
<tr><th></th><th>Replicator</th><th>BNN</th><th>Logit</th></tr>
<tr><td>$F(v)$</td><td colspan="2">$\dfrac{v^2}{2}$</td><td>$v\ln v - v + 1$</td></tr>
<tr><td>$F'^{(-1)}(u)$</td><td colspan="2">$u$</td><td>$e^u$</td></tr>
<tr><td>$\rho$</td><td colspan="2">$(v)_+$, $v \in \mathbb{R}$</td><td>$\rho\left(U(\cdot,\mu_t)\right) = \dfrac{e^{U(\cdot,\mu_t)\eta^{-1}}}{\int_\Omega e^{U(y,\mu_t)\eta^{-1}}\mathrm{d}y}$</td></tr>
<tr><td>$\lambda_t(\mathrm{d}x)$</td><td>$\mu_t(\mathrm{d}x)$</td><td colspan="2">Some $\kappa \in \mathcal{P}$ with a full support on $\Omega$</td></tr>
</table>

### 2.2.2 Specific setting: difference maximization case

We assume that $F$ is strictly increasing. Afterward, we can solve (3) as follows with a Lagrangian multiplier $\eta > 0$:

$$u_x^*(z) = F'^{(-1)}\left(\frac{\left(U(z,\sigma)-U(x,\sigma)\right)_+}{\eta}\right), \tag{4}$$

where $F'$ is the differential of $F$ and $F'^{(-1)}$ is the inverse of $F'$. The utility difference appears on the right-hand side of (4) but with a scaling by $\eta$. Then, the right-hand side of (4) can be understood as a protocol $\rho$ when we substitute $\sigma = \mu_t$, so that it is linked to a maximization principle. Note that the approach by Yoshioka [49] is essentially the formulation with $\eta = 1$. The cost in this case implies that there is no cost if there is no gain by updating the utilities ($U(z,\sigma) \leq U(x,\sigma)$). Conversely, a cost is incurred if agents want to update their actions so that there is some gain ($U(z,\sigma) > U(x,\sigma)$). These observations, along with the cost constraint, imply that searching for better decisions to relatively increase the utility is costly in the difference maximization case.

The multiplier $\eta$ should depend on time $t$ when $\sigma = \mu_t$ in (4) and is therefore rewritten as $\eta_t$. If we find a suitable function $F$, then we can consider the following new EGD instead of (2):

$$\begin{aligned} \frac{\mathrm{d}\mu_t(A)}{\mathrm{d}t} = & \int_{x\in A}\int_{z\in\Omega} F'^{(-1)}\left(\frac{\left(U(x,\mu_t)-U(z,\mu_t)\right)_+}{\eta_t}\right)\lambda_t(\mathrm{d}x)\mu_t(\mathrm{d}z) \\ & -\int_{x\in A}\int_{z\in\Omega} F'^{(-1)}\left(\frac{\left(U(z,\mu_t)-U(x,\mu_t)\right)_+}{\eta_t}\right)\lambda_t(\mathrm{d}z)\mu_t(\mathrm{d}x) \end{aligned} . \tag{5}$$

Note that the multiplier $\eta_t$ depends implicitly on $\mu_t$ because of the optimization formulation (3), and (5) is therefore more strongly nonlinear than (2).

***Remark 1*** For the BNN and replicator models, a suitable $F$ is $F(v) = v^2/2$ ($v \geq 0$).

### 2.2.3 Specific setting: value maximization case

The logit-type protocol $\rho$ has already been suggested to arise from a maximization principle (e.g., Eq. (3) in Lahkar et al. [29]): For each $x \in \Omega$ and $\mu \in \mathcal{P}(\Omega)$,

$$\frac{e^{\eta^{-1}U(\cdot,\mu)}}{\int_\Omega e^{\eta^{-1}U(y,\mu)}\mathrm{d}y} = \underset{u(\cdot)\mathrm{d}x\in\mathcal{P}(\Omega)}{\arg\max}\left\{\underbrace{\int_\Omega U(z,\mu)u(z)\mathrm{d}z}_{\text{Net utility}} - \eta\underbrace{\int_\Omega\left\{u(z)\ln u(z)-\left(u(z)-1\right)\right\}\mathrm{d}z}_{\text{Cost}}\right\}. \tag{6}$$

This is essentially a dual of (3) when $F(v) = v\ln v - v + 1$ ($v > 0$) with $F(0) = 0$ because this $F$ is strongly convex and nonnegative.

As shown later (**Proof of Proposition 2**), for each $\varepsilon > 0$, there exists a unique $\eta > 0$ such that the left-hand side of (6) becomes the maximizer in (3) with equality. Moreover, this $\eta$ is decreasing in $\varepsilon$, showing that paying a larger cost allows for choosing a smaller $\eta$, which gives $\mu$ that is closer to a Nash equilibrium (e.g., Lahkar et al. [29]) (**Section A.1** presents the definition of Nash equilibria). The cost in this case is the relative entropy between $u$ and the uniform distribution on $\Omega$, meaning that there is no cost if we choose ($F(1) = 0$) the uniform distribution, with which agents completely randomly update their actions irrespective of their utilities. Conversely, choosing actions that depend on utilities is costly (needs more effort for better exploration), and updating actions to probably yield a better outcome is more costly in the value maximization case.

Consequently, by explicating the time-dependence of the Lagrangian multiplier $\eta = \eta_t$, we can consider the following new EGD instead of (2):

$$\frac{\mathrm{d}\mu_t(A)}{\mathrm{d}t} = \frac{\int_A e^{\eta_t^{-1}U(y,\mu_t)}\mathrm{d}y}{\int_\Omega e^{\eta_t^{-1}U(y,\mu_t)}\mathrm{d}y} - \mu_t(A) \tag{7}$$

with the protocol

$$\rho(U(z,\mu_t)) = \frac{e^{\eta_t^{-1}U(z,\mu_t)}}{\int_\Omega e^{\eta_t^{-1}U(y,\mu_t)}\mathrm{d}y}, \quad z \in \Omega. \tag{8}$$

Again, $\eta_t$ depends on $\mu_t$ because of the optimization formulation, and (7) is more strongly nonlinear than (2) (and hence than the classical logit dynamic where $\eta$ is simply a constant).

### 2.3 Introducing the value function: forward-looking behavior

In most cases, the utility $U = U(x,\mu_t)$ is determined solely based on the current information; that is, it is given without any future perspectives. Common EGDs assume the following utility [e.g.,13,19,27,28]:

$$U(x,\mu_t) = \int_{y\in\Omega} f_1\left(x, y, \int_{z\in\Omega} z\mu_t(\mathrm{d}z)\right)\mu_t(\mathrm{d}y) + f_2(x) \tag{9}$$

with a trivariate function $f_1$ and a univariate function $f_2$ with $\int_{z\in\Omega} z\mu_t(\mathrm{d}z)$ being the average action; indeed, the right-hand side of (9) is simply an integration of a function by the current probability measure $\mu_t$ of agent actions.

We want to extend the utility $U = U(x,\mu_t)$ so that it depends on some perspective of agents, i.e., prediction in the future $s > t$. This issue is addressed in this paper by replacing $U$ with an expectation conditioned on the current information. More specifically, (2) can be understood as a Fokker–Planck equation (often called the forward Kolmogorov equation) of the following mean-field stochastic differential equation (SDE) that governs a scalar process $X = (X_t)_{t\geq 0}$ in a filtered probability space, which is a time-dependent action of a representative agent [e.g.,49,51]:

$$\mathrm{d}X_t = \mathrm{d}Z_t, \quad t > 0 \tag{10}$$

subject to an initial condition $X_0 \in \Omega$. Here, $Z = (Z_t)_{t\geq 0}$ is a càdlàg (right-continuous and possesses left limits) jump process with the jump measure $\rho(U(z,\mu_t), U(X_t,\mu_t))\lambda_t(\mathrm{d}z)$, and its value is updated at each jump time $\tau$ as follows ("$\tau -$" means the left-limit at time $\tau$):

$$Z_\tau = V_\tau - X_{\tau-} + Z_{\tau-} \quad \text{so that} \quad X_\tau = V_\tau, \tag{11}$$

where $V_\tau$ is generated from the probability measure $c_\rho \rho\left(U\left(z,\mu_{\tau-}\right),U\left(X_{\tau-},\mu_{\tau-}\right)\right)\lambda_{\tau-}(\mathrm{d}z)$, with $c_\rho$ being a normalization constant. Here, recall that $\mu$ represents the distribution of agent actions and is here identified as the distribution of $X$.

With the SDE (10) in mind, at each time $t>0$, one may consider the following (uncontrolled) value function $\Phi_0:\Omega\to\mathbb{R}$, which is a conditional expectation:

$$\Phi_0(x)=\mathbb{E}\left[\delta\int_t^{+\infty}U\left(X_s,\mu_s\right)e^{-\delta(s-t)}\mathrm{d}s\middle|X_t=x\right],\quad x\in\Omega \tag{12}$$

with a discount rate $\delta>0$. This is a forward-looking version of $U$, where the strength of the forward-looking nature can be modulated by $\delta$ such that selecting a smaller value of $\delta$ results in a larger weighting in the future. This is qualitatively the same theory development as Iijima and Oyama [23].

Now, we want to incorporate an optimization mechanism subject to cost constraints into (12) in a tractable way. In this paper, we propose the following version:

$$\Phi(x)=\sup_\rho\mathbb{E}\left[\delta\int_t^{+\infty}U\left(X_s,\mu_t\right)e^{-\delta(s-t)}\mathrm{d}s-\eta\int_s^{+\infty}F\left(\rho_s\right)e^{-\delta(s-t)}\mathrm{d}s\middle|X_t=x\right],\quad x\in\Omega,\quad t>0 \tag{13}$$

subject to the SDE (10), where the second integral in (13) represents the exploration cost and $\rho$ represents the control variable corresponding to the protocol, i.e., jump rate. Note that we use the current distribution $\mu_t$ in the first integral of (13) so that it is a given quantity here. This corresponds to a situation where agents predict utility in the future, the value function, based on the current utility function $U(\cdot,\mu_t)$. The benefit of this formulation is that the right-hand side of (13) is seen as a maximization problem in stochastic control to which the dynamic programming argument applies; thus, the governing equation of the value function $\Phi$ is the static HJB equation (e.g., Chapter 5 in Øksendal and Sulem [39]):

$$\delta\Phi(x)=\delta U\left(x,\mu_t\right)+\sup_{u_x}\left\{\int_\Omega u_x(z)\left(\Phi(z)-\Phi(x)\right)\lambda_t(\mathrm{d}z)-\eta\int_\Omega F\left(u_x(z)\right)\lambda_t(\mathrm{d}z)\right\},\quad x\in\Omega,\quad t>0. \tag{14}$$

This HJB equation (14) still has a free parameter $\eta$, the Lagrangian multiplier, which is determined by the following cost constraint:

$$\int_{x\in\Omega}\int_{z\in\Omega}F\left(u_x^*(z)\right)\lambda_t(\mathrm{d}z)\lambda_t(\mathrm{d}x)\le\varepsilon,\quad t>0, \tag{15}$$

where $u_x^*$, which depends on $\eta$, is the maximizer on the right-hand side of (14). The maximization problem formed by (14)-(15) is solved in exactly the same way as those conducted in **Section 2.2**. We then obtain the counterparts of (4) and (8), where the utility $U(\cdot,\mu_t)$ is replaced by the value function $\Phi(\cdot)$; here, note that $\Phi$ implicitly depends on $\mu_t$.

Now, we explain the ideas for accounting for forward-looking behavior and exploration cost in EGDs and how to combine them based on an HJB equation. Our framework is based on the coupling between the EGD (2), where the utility is replaced by the value function $\Phi$, and the HJB equation (14) along with the constraint (15). The EGD is a time-dependent PIDE, while the HJB equation is a static equation whose coefficients depend on time. In this way, our formulation avoids using a forward-backward system that is typically met in conventional MFGs. Nevertheless, the resulting system is still more complex than classical EGDs because $\Phi$ and $\eta$ depend on $\mu$.

***Remark 2*** A comparison between the static optimization problem (3) and the HJB equation (14) with constraint (15) shows that the maximization mechanism in the latter is analogous to that of the former, where $U$ is effectively replaced by $\Phi$. Moreover, at least formally, the latter approaches the former under the limit $\delta \to +\infty$ if we can expect $\Phi(\cdot) = U(\cdot, \mu_t)$ under this limit due to (14).

## 3. Mathematical analysis

### 3.1 Preparation

We formulate a new EGD that accounts for forward-looking behavior and cost constraints. We first present a general formulation and then separately address the difference maximization case and value maximization case because they are technically different from each other.

In **Section 3**, we use the following notations. The space of all continuous functions on $\Omega$ is denoted by $C(\Omega)$, which is equipped with the usual maximum norm $\|f\|_\infty = \sup_{x\in\Omega}|f(x)|$ for any $f \in C(\Omega)$. The space of all the signed measures on $\Omega$ is denoted by $\mathcal{M}(\Omega)$, which is equipped with the total variation norm $\|\mu\|_{\mathrm{TV}} = \sup_{|f|\le 1}\left|\int_\Omega f(x)\mu(\mathrm{d}x)\right|$ for any $\mu \in \mathcal{M}(\Omega)$. Here, the supremum is taken with respect to all real-valued bounded functions $f$ on $\Omega$. The space of all the signed measures on $\Omega$ whose total variation norm is not greater than 2 is denoted by $\mathcal{M}_2(\Omega) = \left\{\mu \in \mathcal{M}(\Omega);\ \|\mu\|_{\mathrm{TV}} \le 2\right\}$. Note that $\mathcal{P}(\Omega) \subset \mathcal{M}_2(\Omega)$. For simplicity, we omit "$(\Omega)$" from the notations of functional spaces; e.g., $C(\Omega)$ is written simply as $C$ and $\mathcal{P}(\Omega)$ as $\mathcal{P}$.

Throughout this paper, we employ the following assumption about utility, some of which **((A2)-(A4)**) have bene commonly employed in the study of EGDs.

***Assumption 1***

*Utility* $U : \Omega \times \mathcal{M}_2 \to \mathbb{R}$ *satisfies the following conditions:*

***(A1) Profile.*** *For all* $\mu \in \mathcal{M}_2$, $U(\cdot, \mu)$ *has a finite number of isolated global maxima, and for each maximizing point* $x'$, $U(\cdot, \mu)$ *is twice continuously differentiable and there is a constant* $c > 0$ *with* $U(x', \mu) - U(x, \mu) \approx c(x - x')^2$ *in a neighborhood of* $x'$. *Moreover, there is a constant* $\Delta U > 0$ *with*

$$\inf_{\mu\in\mathcal{M}_2}\left\{\sup_{x\in\Omega} U(x,\mu) - \inf_{x\in\Omega} U(x,\mu)\right\} \ge \Delta U. \tag{16}$$

***(A2) Boundedness.***

$$0 \le U(x, \mu) \le \bar{U} \tag{17}$$

*with a constant* $\bar{U} > 0$ *for all* $(x, \mu) \in \Omega \times \mathcal{M}_2$,

***(A3) Lipschitz continuity 1.*** $\left|U(x,\mu)-U(x,\sigma)\right| \le L_U \left\|\mu-\sigma\right\|_{\mathrm{TV}}$ *(18)*

*with a constant* $L_U > 0$ *for all* $(x,\mu,\sigma) \in \Omega \times \mathcal{M}_2 \times \mathcal{M}_2$,

***(A4) Lipschitz continuity 2.*** $\left|U(x,\mu)-U(y,\mu)\right| \le \mathrm{Lip}_U \left|x-y\right|$ *(19)*

*with a constant* $\mathrm{Lip}_U > 0$ *for all* $(x,y,\mu) \in \Omega \times \Omega \times \mathcal{M}_2$, *and*

***(A5) Entropic condition.*** *Fix* $\varepsilon > 0$. *For all* $\mu \in \mathcal{M}_2$, *any solution* $m = m(\mu) > 0$ *to the following equation is bounded below by a constant* $\underline{\eta}_{\varepsilon} > 0$ *independent of* $\mu$ *:*

$$\varepsilon = \int_{z\in\Omega} \left( \frac{\exp\left(m(\mu)^{-1} U(z,\mu)\right)}{\int_{\Omega} \exp\left(m(\mu)^{-1} U(y,\mu)\right) \mathrm{d}y} \ln \frac{\exp\left(m(\mu)^{-1} U(z,\mu)\right)}{\int_{\Omega} \exp\left(m(\mu)^{-1} U(y,\mu)\right) \mathrm{d}y} \right) \mathrm{d}z \,. \tag{20}$$

We need to extend the domain of utilities from $\Omega \times \mathcal{P}$ to $\Omega \times \mathcal{M}_2$ in order to obtain theoretical results about our EGD. A trivial utility that satisfies **Assumption 1** is $U(x,\cdot) = f_2(x)$ with a suitable nonconstant function $f_2$, which is included in (9). One may also consider (9) with $\mu$ replaced by its positive part $\mu^+$; taking the positive part is not harmful in practice because we eventually consider $\mu \in \mathcal{P}$.

The assumption **(A1)** means that any utility is not a constant function and is not flat around its maxima. This is technically necessary for our unique existence result of solutions to forward-looking EGDs, particularly to ensure the unique existence of a solution that satisfies the cost constraint (**Proof of Proposition 2**). The assumptions **(A2)-(A3)** ensure that there is neither a blow-up nor an extremely rapid variation in the utility by changing the distribution of agent actions. The nonnegativity of utilities is not restrictive; if a utility is bounded but possibly becomes negative, then we can add a sufficiently large positive constant to it without any loss of generality. The assumption **(A4)** means that the utility is a continuous function for each distribution of agent actions and is not allowed to vary rapidly in $x$. The final assumption **(A5)** seems to be the most complicated, but the right-hand side of (20) is essentially the minus of entropy, and this condition means that the probability density proportional to $\exp\left(m(\mu)^{-1} U(z,\mu)\right)$ is moderately different from the uniform distribution.

We consider that **(A2)-(A4)** are not restrictive because they are satisfied in many examples [e.g., 13,19,27,28]. However, there exists an example where it fails, such as the continuity in $x$ [e.g., 2]. In such a case, one may approximate the discontinuity with a large Lipschitz constant $\mathrm{Lip}_U$. In contrast, **(A1) and (A5)** are used for analyzing the logit model and are more restrictive and specific to our EGDs. Despite this theoretical restriction, we demonstrate that EGDs can be numerically investigated by a simple scheme.

### 3.2 Problem formulation

For each fixed constant $\underline{\eta} > 0$ and $\varepsilon > 0$, the forward-looking EGD is formulated as follows: Given an initial condition $\mu_0 \in \mathcal{P}(\Omega)$, the solution $\mu = (\mu_t)_{t \geq 0}$ is governed by

$$\begin{aligned}\frac{\mathrm{d}\mu_t(A)}{\mathrm{d}t} &= \int_{x\in A}\int_{z\in\Omega} \rho\left(\frac{\Phi(z,\mu_t)}{\eta(\mu_t)}, \frac{\Phi(x,\mu_t)}{\eta(\mu_t)}\right)\lambda_t(\mathrm{d}x)\mu_t(\mathrm{d}z) \\ &\quad - \int_{x\in A}\int_{z\in\Omega} \rho\left(\frac{\Phi(x,\mu_t)}{\eta(\mu_t)}, \frac{\Phi(z,\mu_t)}{\eta(\mu_t)}\right)\lambda_t(\mathrm{d}z)\mu_t(\mathrm{d}x)\end{aligned} \tag{21}$$

for any $t > 0$ and $A \subset \mathcal{B}$. At each $t > 0$, the couple $\left(\eta(\mu_t), \Phi(\cdot,\mu_t)\right) \in \mathbb{R}_+ \times C$ is governed by the system

$$\delta\Phi(x,\mu_t) = \delta U(x,\mu_t) + \left\{\begin{matrix}\int_\Omega u_x^*(z)\left(\Phi(z,\mu_t) - \Phi(x,\mu_t)\right)\lambda_t(\mathrm{d}z) \\ -\eta(\mu_t)\int_\Omega F\left(u_x^*(z)\right)\lambda_t(\mathrm{d}z)\end{matrix}\right\},\quad x \in \Omega \tag{22}$$

with

$$u_x^*(\cdot) = \underset{u_x(\cdot)\geq 0\ \left(\int_\Omega u_x(z)\mathrm{d}z = 1\right)}{\arg\max} \left\{\int_\Omega u_x(z)\left(\Phi(z,\mu_t) - \Phi(x,\mu_t)\right)\lambda_t(\mathrm{d}z) - \eta(\mu_t)\int_\Omega F\left(u_x(z)\right)\lambda_t(\mathrm{d}z)\right\} \tag{23}$$

s.t.

$$\int_{x\in\Omega}\int_{z\in\Omega} F\left(u_x^*(z)\right)\lambda_t(\mathrm{d}z)\lambda_t(\mathrm{d}x) + \chi R\left(\eta(\mu_t)\right) \leq \varepsilon. \tag{24}$$

Here, the time derivative of $\mu_t$ is understood in the following sense:

$$\lim_{h\to 0}\left\|\frac{\mathrm{d}\mu_t}{\mathrm{d}t} - \frac{\mu_{t+h} - \mu_t}{h}\right\|_{\mathrm{TV}} = 0. \tag{25}$$

The term $\chi R\left(\eta(\mu_t)\right)$ consists of a small positive parameter $\chi > 0$ and a decreasing function $R(\cdot) : (0, +\infty) \to \mathbb{R}_+$ to be introduced to regularize the dynamics.

***Remark 2*** The condition "$\left(\int_\Omega u_x(z)\mathrm{d}z = 1\right)$" is activated only for the logit model (**Section 3.4**).

***Remark 3*** The term $R\left(\eta(\mu_t)\right)$ in (24) is introduced only for the BNN and replicator models (**Section 3.3**). This term with a specific functional form plays a role in well-posing our EGD, with which we can obtain a strict lower bound of $\eta(\mu_t)$, which is necessary to prove the Lipschitz continuity of the term $\eta^{-1}\Phi$ with respect to $\mu_t$. Our computational results in **Section 4** suggest how this term works.

***Remark 4*** We explicitly represent the dependence of the value function $\Phi$ and optimal Lagrangian multiplier $\eta$ on the distribution $\mu_t$.

### 3.3 Difference maximization case: BNN and replicator models

In this case, we reformulate the system (21)-(24) as follows, where we set $\lambda_t(\mathrm{d}x) = w\kappa(\mathrm{d}x) + (1-w)\mu_t(\mathrm{d}x)$ with some $\kappa \in \mathcal{P}(\Omega)$ fully supported on $\Omega$ and a constant $w \in [0,1]$ to address BNN and replicator models in a unified way: For each $\varepsilon > 0$ and $w \in [0,1]$, given an initial condition $\mu_0 \in \mathcal{P}(\Omega)$, the solution $\mu = (\mu_t)_{t \geq 0}$ is governed by the EGD

$$\begin{aligned} \frac{\mathrm{d}\mu_t(A)}{\mathrm{d}t} = & \int_{x \in A} \int_{z \in \Omega} \left( \frac{\Phi(x,\mu_t) - \Phi(z,\mu_t)}{\eta(\mu_t)} \right)_+ \lambda_t(\mathrm{d}x) \mu_t(\mathrm{d}z) \\ & - \int_{x \in A} \int_{z \in \Omega} \left( \frac{\Phi(z,\mu_t) - \Phi(x,\mu_t)}{\eta(\mu_t)} \right)_+ \lambda_t(\mathrm{d}z) \mu_t(\mathrm{d}x) \end{aligned} \tag{26}$$

for any $t > 0$ and $A \subset \mathcal{B}$, where $(\eta(\mu_t), \Phi(\cdot,\mu_t)) \in \mathbb{R}_+ \times C$ is governed by the system

$$\delta\Phi(x,\mu_t) = \delta U(x,\mu_t) + \left\{ \begin{array}{l} \int_\Omega u_x^*(z) \left( \Phi(z,\mu_t) - \Phi(x,\mu_t) \right) \lambda_t(\mathrm{d}z) \\ -\eta(\mu_t) \int_\Omega \frac{1}{2} \left( u_x^*(z) \right)^2 \lambda(\mathrm{d}z) \end{array} \right\}, \quad x \in \Omega \tag{27}$$

with

$$u_x^*(\cdot) = \underset{u_x(\cdot) \geq 0}{\arg\max} \left\{ \int_\Omega u_x(z) \left( \Phi(z,\mu_t) - \Phi(x,\mu_t) \right) \lambda_t(\mathrm{d}z) - \eta(\mu_t) \int_\Omega \frac{1}{2} \left( u_x(z) \right)^2 \lambda_t(\mathrm{d}z) \right\} \tag{28}$$

s.t.

$$\int_{x \in \Omega} \int_{z \in \Omega} \left( \frac{1}{2} \left( u_x^*(z) \right)^2 \right) \lambda_t(\mathrm{d}z) \lambda_t(\mathrm{d}x) + \frac{\chi}{\eta(\mu_t)^{2+\xi}} = \varepsilon, \tag{29}$$

where $\xi \geq 0$ and $\chi > 0$ are constants.

The following proposition states the unique existence of the forward-looking EGD (26), which covers both BNN and replicator models. Here, we prove the case $\xi = 0$ that suffices for theoretical purposes here and investigate cases $\xi > 0$ in **Section 4** to compare regularization methods.

***Proposition 1***

*Assume that $\delta > 0$ is sufficiently large and $\xi = 0$. Then, the forward-looking EGD (26) admits a unique solution $\mu_t \in \mathcal{P}$ ($t \geq 0$). At each $t > 0$, it follows that $0 \leq \Phi(x,\mu_t) \leq \bar{U}$ ($x \in \Omega$). Moreover, $\underline{\eta} \leq \eta(\sigma) \leq \bar{\eta}$ with some constants $\underline{\eta}, \bar{\eta} > 0$ for all $\sigma \in \mathcal{P}$.*

The strategy to prove **Proposition 1** is as follows. We apply Theorem 1 in Cheung [12] to the forward-looking EGD (26), and we show that it suffices to prove that the quantity $\eta(\sigma)^{-1}\Phi(\cdot,\sigma)$ at each $\Omega$ is bounded and Lipschitz continuous with respect to $\sigma \in \mathcal{P}$. We address this issue by showing that $\Phi$ is uniformly bounded irrespective of $(x,\sigma)$ and is Lipschitz continuous with respect to $\sigma \in \mathcal{P}$. We also show that $\eta(\sigma)$, owing to the regularization in (29), is strictly bounded (it is also strictly positive) and Lipschitz continuous with respect to $\sigma \in \mathcal{P}$ if it exists. For each $\sigma \in \mathcal{P}$, the unique existence of the couple $\left(\eta(\sigma),\Phi(\cdot,\sigma)\right)\in\mathbb{R}_{+}\times C$ is proven on the basis of the Banach fixed-point theorem, which completes the proof due to the strict boundedness of them.

In the present case, we consider the cost constraint (29) as an equality but not an inequality because it admits a nonstandard form in view of convex optimization since it contains the Lagrangian multiplier. The regularization term in (29) plays a role here in obtaining a positive lower bound of $\eta$ that does not depend on $\sigma$. This term can be understood as a fixed cost for updating decisions. Investigations of the BNN and replicator models without regularization seem to need another method that does not rely on the strict boundedness of $\eta$. Note that **(A1) and (A5)** are not used in **Proof of Proposition 1**.

### 3.4 Value maximization case: Logit model

In this case, we reformulate the system (21)-(24) as follows: For each $\varepsilon>0$, given an initial condition $\mu_0\in\mathcal{P}(\Omega)$, the solution $\mu=\left(\mu_t\right)_{t\geq 0}$ is governed by the EGD

$$\frac{\mathrm{d}\mu_t(A)}{\mathrm{d}t}=\frac{\int_A\exp\left(\dfrac{\Phi(y,\mu_t)}{\eta(\mu_t)}\right)\mathrm{d}y}{\int_\Omega\exp\left(\dfrac{\Phi(y,\mu_t)}{\eta(\mu_t)}\right)\mathrm{d}y}-\mu_t(A) \tag{30}$$

for any $t>0$ and $A\subset\mathcal{B}$, where $\left(\eta(\mu_t),\Phi_t(\cdot,\mu_t)\right)\in\mathbb{R}_{+}\times C$ is governed by the system

$$\delta\Phi(x,\mu_t)=\delta U(x,\mu_t)+\left\{\begin{array}{l}\int_\Omega u_x^*(z)\left(\Phi(z,\mu_t)-\Phi(x,\mu_t)\right)\mathrm{d}z\\-\eta(\mu_t)\int_\Omega\left\{u_x^*(z)\ln\left(u_x(z)\right)-u_x(z)+1\right\}\mathrm{d}z\end{array}\right\},\quad x\in\Omega \tag{31}$$

with

$$u_x^*(\cdot)=\underset{u_x(\cdot)\geq 0,\ \int_\Omega u_x(z)\mathrm{d}z=1}{\arg\max}\left\{\begin{array}{l}\int_\Omega u_x(z)\left(\Phi(z,\mu_t)-\Phi(x,\mu_t)\right)\mathrm{d}z\\-\eta(\mu_t)\int_\Omega\left\{u_x(z)\ln\left(u_x(z)\right)-u_x(z)+1\right\}\mathrm{d}z\end{array}\right\} \tag{32}$$

s.t.

$$\int_{x\in\Omega}\int_{z\in\Omega}\left\{u_x^*(z)\ln\left(u_x^*(z)\right)-u_x^*(z)+1\right\}\mathrm{d}z\mathrm{d}x\leq\varepsilon\ . \tag{33}$$

The next proposition states a unique existence result of the forward-looking EGD (30) of the logit type.

***Proposition 2***

*Assume that $\delta>0$ is sufficiently large. Then, the forward-looking EGD (30) admits a unique solution $\mu_t\in\mathcal{P}$ ($t\geq 0$). At each $t>0$, it follows that $0\leq\Phi(x,\mu_t)\leq\bar{U}$ ($x\in\Omega$). Moreover, $\underline{\eta}\leq\eta(\sigma)\leq\bar{\eta}$ with some constants $\underline{\eta},\bar{\eta}>0$ for all $\sigma\in\mathcal{P}$.*

The strategy to prove **Proposition 2** is qualitatively the same as **Proposition 1**, but they are technically different at several points. The greatest difference is that the HJB equation (31) is solvable in a closed form for **Proposition 2**. This remarkable property facilitates the analysis of the cost constraint (33), which is actually more complicated than the previous one (29). Moreover, no regularization is necessary in (33). Instead, we need structural conditions **(A1) and (A5)**.

### 3.5 Some remarks

A technical difference between the difference maximization case and the value maximization case is that the HJB equation is not explicitly solvable in the former:

$$\Phi(x) = U(x,\mu_t) + \frac{1}{2\delta\eta(\mu_t)}\int_\Omega \left(\Phi(z)-\Phi(x)\right)_+^2 \lambda_t(\mathrm{d}z),\quad x\in\Omega, \tag{34}$$

while it is possible in the latter:

$$\Phi(x) = \frac{\delta}{\delta+1}U(x,\mu_t) + \frac{\eta(\mu_t)}{\delta}\ln\left(\int_\Omega \exp\left(\frac{\delta}{\delta+1}\frac{U(z,\mu_t)}{\eta(\mu_t)}\right)\mathrm{d}z\right),\quad x\in\Omega. \tag{35}$$

In this case, the value function $\Phi$ is a convex combination between the utility and a utility-based constant. Such a simple structure of $\Phi$ may not exist for the difference maximization case.

We formulated and mathematically analyzed forward-looking EGDs for BNN, replicator, and logit models. A difference between these models is that $\eta$ directly affects stationary states if they exist for the logit model but not at least explicitly for the BNN and replicator models. This is because of the following equations for a stationary state:

$$\int_{x\in A}\int_{z\in\Omega}\left(\Phi(x,\mu_t)-\Phi(z,\mu_t)\right)\lambda_t(\mathrm{d}x)\mu_t(\mathrm{d}z) = 0 \tag{36}$$

for the BNN and replicator models and

$$\mu_t(A) = \frac{\int_A \exp\left(\frac{\Phi(y,\mu_t)}{\eta(\mu_t)}\right)\mathrm{d}y}{\int_\Omega \exp\left(\frac{\Phi(y,\mu_t)}{\eta(\mu_t)}\right)\mathrm{d}y} \tag{37}$$

for the logit model. In (36), $\eta$ does not appear but is implicitly linked to $\Phi$. In contrast, in (36), $\eta$ explicitly appears on the right-hand side of (37). Inspecting the formulae (34) and (35) implies that we formally have $\Phi \approx U$ at $\delta \to +\infty$ because the second terms in these equations are strictly bounded and vanish at $\delta \to +\infty$. Our formulation therefore consistently covers cases with and without mean-field effects. For this myopic case $\delta \to +\infty$, (36) formally reduces to

$$\int_{x\in A}\int_{z\in\Omega}\left(U(x,\mu_t)-U(z,\mu_t)\right)\lambda_t(\mathrm{d}x)\mu_t(\mathrm{d}z) = 0, \tag{38}$$

which is free from $\eta$. Note that equilibria possibly depend on initial conditions with a replicator protocol because the support of an initial condition is preserved for all $t > 0$ (e.g., Theorem 13 in Mendoza-Palacios and Hernández-Lerma [36]). In contrast, under the same limit, (37) reduces to

$$\mu_t(A) = \frac{\int_A \exp\left(\dfrac{U(y,\mu_t)}{\eta(\mu_t)}\right)\mathrm{d}y}{\int_\Omega \exp\left(\dfrac{U(y,\mu_t)}{\eta(\mu_t)}\right)\mathrm{d}y}, \tag{39}$$

which still depends on $\eta$, and any equilibrium has full support on $\Omega$ due to the right-hand side of (39). These observations imply that the optimal multiplier $\eta$ affects mainly the convergence speed for the BNN and replicator models but affects the equilibria for the logit model. Owing to the dependence on $\eta$, an equilibrium in the logit model is different from the Nash equilibria but is possibly their approximation for small $\eta$.

The difference maximization case and value maximization case also differ in terms of the control variables because the former optimizes $\rho$ in the space of nonnegative functions but the latter in the space of probability measures, and the latter space is strictly smaller. This is due to following the formulation of previous studies, and a completely different EGD may arise for the value maximization case if one optimizes the protocol in the space of nonnegative functions with the cost constraint being unchanged.

A marked difference between the classical and our MFG formulations is that the former often needs a monotonicity condition (concavity of utility $U(\cdot,\mu_t)$ with respect to $\mu_t$ in conjunction with some structural condition about $\rho$), whereas our formulation does not; in classical MFGs, monotonicity has been employed to ensure the uniqueness of solutions and is essential for investigating the forward-backward system [10,17]. In contrast, our proofs rely on the fixed-point argument for the HJB part and Lipschitz continuity of utility. The cost constraint has not been included in classical MFG models, which is also a difference from our model.

This paper focuses on the one-dimensional action space, while its multidimensional extension (e.g., Chapter 5 in Mendoza-Palacios and Hernández-Lerma [36]) is also straightforward if the domain is compact and has a simple structure such as a multidimensional cube, as computationally investigated in **Section 4**. We consider that extensions of the proposed EGDs to those played by heterogeneous agents [49] are also possible; however, interactions among agents, which are encoded in utility, protocol, discount, and/or cost constraints, would need care so that the unique existence of solutions is not broken. The conditions to guarantee it would be more severe than in the proposed case. These extensions are beyond the scope of this paper but will be addressed elsewhere.

### 3.6 Numerical discretization

#### 3.6.1 Setting

In this subsection, we explain how to computationally implement forward-looking EGD. We use a finite difference method based on a fixed time increment $\Delta t > 0$ and spatial increment $\Delta x = 1/N$, as in existing studies (e.g., Chapter 5 in Mendoza-Palacios and Hernández-Lerma [36]; Yoshioka and Tsujimura [50]), where $N \in \mathbb{N}$ with $N \geq 2$ is the degree of freedom in the $x$ direction. We prepare cells $\Omega_i = \left[ (i-1)\Delta x, i\Delta x \right)$ for $i = 1,2,...,N-1$ and $\Omega_N = \left[ 1-\Delta x, 1 \right]$. Each quantity discretized at time step $k\Delta t$ ($k = 0,1,2,...$) and in $\Omega_i$ ($i = 1,2,...,N$) is denoted by using subscript $i,k$, e.g., $\mu_{i,k}$ and $\Phi_{i,k}$. The initial condition $\mu_0$ is discretized in $\Omega_i$ as $\mu_{0,i} = \int_{\Omega_i} \mu_0 (\mathrm{d}x)$ or its approximation if necessary. We explain the discretization for the difference maximization case (26) of the replicator type and value maximization case (30). Discretization for the BNN model is essentially the same as that for replicator one and is therefore omitted.

### 3.6.2 Replicator model

Assume that we already have $\mu_{i,k-1}$ ($i = 1,2,...,N$) for some $k \in \mathbb{N}$. Each $\mu_{i,k}$ ($i = 1,2,...,N$) is computed by applying a fully explicit Euler discretization:

$$\mu_{i,k} = \mu_{i,k-1} + \Delta t \sum_{j=1}^{N} \left( \frac{\Phi_{i,k-1} - \Phi_{j,k-1}}{\eta_{k-1}} \right)_+ \mu_{j,k-1}\mu_{i,k-1} - \Delta t \sum_{j=1}^{N} \left( \frac{\Phi_{j,k-1} - \Phi_{i,k-1}}{\eta_{k-1}} \right)_+ \mu_{j,k-1}\mu_{i,k-1} . \tag{40}$$

The quantities $\Phi_{i,k-1}$ ($i = 1,2,...,N$) and $\eta_{k-1}$ are found in the following system on the basis of (34):

$$\Phi_{i,k-1} = U_{i,k-1} + \frac{1}{2\eta_{k-1}\delta} \sum_{j=1}^{N} \left( \Phi_{j,k-1} - \Phi_{i,k-1} \right)_+^2 \mu_{j,k-1} , \quad i = 1,2,...,N \tag{41}$$

and

$$\varepsilon = \frac{1}{2\left(\eta_{k-1}\right)^2} \sum_{i,j=1}^{N} \left( \Phi_{j,k-1} - \Phi_{i,k-1} \right)_+^2 \mu_{i,k-1}\mu_{j,k-1} + \frac{\chi}{\left(\eta_{k-1}\right)^{2+\xi}} , \tag{42}$$

where $U_{i,k-1}$ ($i = 1,2,...,N$) is the discretized utility that is bounded between 0 and $\bar{U}$, and is computed using $\mu_{i,k-1}$ ($i = 1,2,...,N$) and classical midpoint rules evaluated at centers of cells. Our discretization method is therefore an explicit method coupled with a nonlinear algebraic system (41)-(42).

The system (41)-(42) is then iteratively solved by using **Algorithm 1**. The equilibrium point of the system (41)-(42) exists uniquely according to **Proof of Proposition 1** if $\delta > 0$ is large (it suffices to consider cases where $\mu$ is discrete).

### 3.6.3 Logit model

Assume that we already have $\mu_{i,k-1}$ ($i = 1,2,...,N$) for some $k \in \mathbb{N}$. Each $\mu_{i,k}$ ($i = 1,2,...,N$) is computed by applying a fully explicit Euler discretization:

$$\mu_{i,k} = \left(1 - \Delta t\right) \mu_{i,k-1} + \Delta t \frac{\exp\left( \dfrac{\Phi_{i,k-1}}{\eta_{k-1}} \right)}{\sum_{j=1}^{N} \exp\left( \dfrac{\Phi_{j,k-1}}{\eta_{k-1}} \right) \Delta x} . \tag{43}$$

The quantities $\Phi_{i,k-1}$ ($i = 1,2,...,N$) and $\eta_{k-1}$ are found from the following system, where $\Phi_{i,k-1}$ is found explicitly from (35):

$$\Phi_{i,k-1} = U_{i,k-1} + \frac{\eta_{k-1}}{\delta} \ln\left( \sum_{j=1}^{N} \exp\left( \frac{\delta}{\delta+1} \frac{U_{j,k-1}}{\eta_{k-1}} \right) \Delta x \right), \quad i = 1,2,\ldots,N \tag{44}$$

and

$$\varepsilon = \sum_{i=1}^{N} \left\{ \frac{\exp\left( \frac{\delta}{\delta+1} \frac{U_{i,k-1}}{\eta_{k-1}} \right)}{\sum_{j=1}^{N} \exp\left( \frac{\delta}{\delta+1} \frac{U_{j,k-1}}{\eta_{k-1}} \right) \Delta x} \ln\left( \frac{\exp\left( \frac{\delta}{\delta+1} \frac{U_{i,k-1}}{\eta_{k-1}} \right)}{\sum_{j=1}^{N} \exp\left( \frac{\delta}{\delta+1} \frac{U_{j,k-1}}{\eta_{k-1}} \right) \Delta x} \right) \right\} \Delta x \, . \tag{45}$$

The equation (45) is solved for $\eta_{k-1}$ by using **Algorithm 2**, where $g_{k-1}^{(n)}$ is expressed as

$$g_{k-1}^{(n)} = \frac{\delta}{\delta+1} \frac{1}{\varepsilon + \ln\left( \sum_{j=1}^{N} \exp\left( \frac{\delta}{\delta+1} \frac{U_{j,k-1}}{\eta_{k-1}^{(n)}} \right) \Delta x \right)} \frac{\sum_{i=1}^{N} U_{i,k-1} \exp\left( \frac{\delta}{\delta+1} \frac{U_{i,k-1}}{\eta_{k-1}^{(n)}} \right) \Delta x}{\sum_{j=1}^{N} \exp\left( \frac{\delta}{\delta+1} \frac{U_{j,k-1}}{\eta_{k-1}^{(n)}} \right) \Delta x} \, . \tag{46}$$

**Algorithm 2** is designed to handle the implicit equation (45) in an iterative way. The equilibrium point of (45) exists uniquely according to **Proof of Proposition 2** if $\delta > 0$ is sufficiently large (it suffices to consider a case where $\mu$ is discrete).

***Remark 5*** In both replicator and logit models, choosing a sufficiently small $\Delta t > 0$ yields nonnegative discrete probability densities at each time step. Both the upper and lower bounds of the discrete value functions hold true as well when discrete $\mu$ cases are considered as in **Proofs of Propositions 1-2**. Additionally, the unique existence of numerical solutions is guaranteed for a sufficiently large $\delta > 0$, which again follows in an analogous way as those in these proofs.

***Remark 6*** Under the setting of **Section 4**, at each time step, **Algorithms 1 and 2** are terminated at most approximately 1000 iterations. Although the exploration of more efficient algorithms that outperform them is beyond the scope of this paper, one may alternatively use another method such as a Newton-like method.

***Algorithm 1***

***Step 0.*** *Set the error threshold* $\mathrm{Er}\left(=10^{-10}\right)$ *and relaxation parameter* $r\left(=0.05\right)$.

***Step 1.*** *Set an initial guess* $\Phi_{i,k-1}^{(0)}=U_{i,k-1}$ *at all grid points and* $\eta_{k-1}^{(0)}>0$, *and a small parameter* $r\in\left(0,1\right]$. *Set* $n=0$.

***Step 2.*** *Compute* $\Phi_{i,k-1}^{(n+1)}=U_{i,k-1}+\frac{1}{2\eta_{k-1}^{(n)}\delta}\sum_{j=1}^{N}\left(\Phi_{j,k-1}^{(n)}-\Phi_{i,k-1}^{(n)}\right)_{+}^{2}\mu_{j,k-1}$ *at all grid points.*

***Step 3.*** *Compute* $\eta_{k-1}^{(n+1)}=r\eta_{k-1}^{(n)}+\left(1-r\right)\sqrt{\frac{1}{2\varepsilon}\sum_{i,j=1}^{N}\left(\Phi_{j,k-1}^{(n)}-\Phi_{i,k-1}^{(n)}\right)_{+}^{2}\mu_{i,k-1}\mu_{j,k-1}+\frac{\chi}{\varepsilon\left(\eta_{k-1}\right)^{\xi}}}$.

***Step 4.*** *Compute the error* $\mathrm{Er}^{(n)}=\max\left\{\max_{1\le i\le N}\left\{\left|\Phi_{i,k-1}^{(n+1)}-\Phi_{i,k-1}^{(n)}\right|\right\},\left|\eta_{k-1}^{(n+1)}-\eta_{k-1}^{(n)}\right|\right\}$.

***Step 4.*** *If* $\mathrm{Er}^{(n)}\le\mathrm{Er}$, *then output* $\Phi_{i,k-1}^{(n+1)}$ *at all grid points and* $\eta_{k-1}^{(n+1)}$ *as the numerical solution to the system (41)-(42) and terminate the algorithm. If* $\mathrm{Er}^{(n)}>\mathrm{Er}$, *then set* $n\to n+1$ *and go to* ***Step 2****.*

***Algorithm 2***

***Step 0.*** *Set the error threshold* $\mathrm{Er}\left(=10^{-10}\right)$ *and relaxation parameter* $r\left(=0.05\right)$.

***Step 1.*** *Set* $\eta_{k-1}^{(0)}>0$ *and a small parameter* $\Delta s>0$. *Set* $n=0$.

***Step 2.*** *Compute* $\eta_{k-1}^{(n+1)}=r\eta_{k-1}^{(n)}+\left(1-r\right)g_{k-1}^{(n)}$.

***Step 3.*** *Compute the error* $\mathrm{Er}^{(n)}=\left|\eta_{k-1}^{(n+1)}-\eta_{k-1}^{(n)}\right|$.

***Step 4.*** *If* $\mathrm{Er}^{(n)}\le\mathrm{Er}$, *then output* $\eta_{k-1}^{(n+1)}$ *as the numerical solution to* (45) *and terminate the algorithm. If* $\mathrm{Er}^{(n)}>\mathrm{Er}$, *then set* $n\to n+1$ *and go to* ***Step 2****.*

## 4. Demonstrative computation

We apply forward-looking EGDs to potential games whose Nash equilibria are known. We also apply them to a two-dimensional problem.

### 4.1 Potential games

#### 4.1.1 Problem setting

We apply forward-looking EGDs to potential games to investigate their performance. We consider the following two utilities: For any $\mu \in \mathcal{P}$, we consider a quadratic utility

$$U_1(x,\mu) = \int_\Omega (x-y)^2 \mu(\mathrm{d}y) \tag{47}$$

and an affine utility

$$U_2(x,\mu) = \left( f\left( \int_\Omega y\mu(\mathrm{d}y) \right) - c \right) x, \tag{48}$$

where $c > 0$ is a constant and $f : \mathbb{R}_+ \to \mathbb{R}_+$ is a strictly positive and decreasing function, which are given here by $c = 2$ and $f(v) = 1/\sqrt{|v| + \varepsilon'}$ ($v \in \mathbb{R}$), where $0 < \varepsilon' << 1$ is to avoid divisions by zero but is neglected in our computation because $\mu$ is approximated at cell centers; hence, any divisions by zero do not occur in the finite difference method. Without any loss of generality, we add a positive constant to (47)-(48), e.g., 1.5 for the latter, to computationally make the utilities positive.

The first utility $U_1$ represents the situation where agent actions should have the greatest variance in $\Omega$ and is considered an idealized problem for testing models. Elementary computations show that a unique Nash equilibrium exists for the quadratic case, which is given by the sum of the Dirac deltas concentrated at the boundary points $x = 0,1$ with equal weight $1/2$. In contrast, the second case is a simple model representing common pool resource management, where $x$ is the harvesting intensity of a resource and the first term $f\left( \int_\Omega y\mu(\mathrm{d}y) \right) x$ of $U_2$ is the utility by harvesting, which decreases as the average harvesting intensity increases. The second term $-cx$ is the harvesting cost. The two terms are assumed to be proportional to $x$, which is not only for simplicity but also because of the following interesting property. It was shown that the Nash equilibria of $U_2$ are not unique and are given by the set of all $\mu \in \mathcal{P}$ such that $\int_\Omega y\mu(\mathrm{d}y) = c^{-2} = 1/4 - \varepsilon'$ (Proposition 4 in Cheung and Lahkar [13]); hence, solutions to EGDs with different initial conditions or protocols may approach different equilibria. Indeed, infinitely many $\mu \in \mathcal{P}$ satisfy this expectation constraint. Therefore, this utility is an interesting case despite its simplicity.

We set $\Delta t = 0.005$ and $\Delta x = N^{-1} = 250$. With this resolution, we can visually observe differences among protocols with reasonably converged numerical solutions (**Section A.4**). We judge a numerical solution to be stationary if the difference between the probability densities $p$, computed by $\mu_{i,k} / \Delta x$ at each time instance and cell, at successive time steps decreases from $10^{-10}$. The initial condition $\mu_0$ is given by a uniform distribution on $\Omega$. For the BNN model, we assume a uniform distribution on $\Omega$ for $\lambda$.

#### 4.1.2 Computational results and discussion

***Behavior of*** $\eta$

We examine the behavior of numerical solutions by focusing on the optimal multiplier $\eta$ for the logit model that does not have any regularization terms in the cost constraint. Here, we set $\delta = 1$ and $\varepsilon = 0.375$ unless otherwise specified. **Figure 1** compares the computed time history of $\eta$ for different values of the cost $\varepsilon$ against the two utilities $U_1$ and $U_2$. **Figure 2** compares the computed time history of the probability density $p$ for $U_1$ and $U_2$. For the first utility $U_1$, **Figure 1(a)** shows that the computed $\eta$ is constant in time for all the cases examined, demonstrating that there are utilities that computationally preserve the exploration cost for the specific initial condition, with which fixing $\eta$ is optimal. Although not presented here, we also observed a similar behavior of the utility $U = -U_1$, but did not observe it for non-uniform initial conditions such as $\mu_0(0, x) = x^2$; see **Section A.5** for the theoretical justification of this phenomenon. For $U_1$, the numerical solutions become closer to being bimodal, and their modes are sharper for larger $\varepsilon$ values with which $\eta$ is smaller; hence, the corresponding numerical solutions are considered better at approximating the Nash equilibrium.

For the second utility $U_2$, **Figure 1(b)** shows that the computed $\eta$ is decreasing in time, suggesting that the optimal strategy for tuning $\eta$ is to gradually decrease it toward some positive level. The decrease in $\eta$ overtime implies that agents use logit functions (the first term in (30)) that depend sharply on utilities as $\mu$ approaches equilibrium. All the computed equilibria have qualitatively the same profiles that are decreasing and convex in $x$ but have different averages of agent actions, as shown in **Table 1**. These average values become closer to the exact value of 0.25 for larger $\varepsilon$. The convergence rate of the average toward the exact one is estimated to be better than the first-order accuracy.

**Figure 3** compares the computed time history of $\eta$ for different values of the discount rate $\delta$ against the two utilities $U_1$ and $U_2$. For both utilities, increasing $\delta$ yields larger values of $\eta$, implying that more myopic agents choose more blurred logit functions that are less sensitive to utilities. Introducing perspectives into agent actions therefore enables us to model EGDs with different sensitivities against utility values. We also investigate this point in a two-dimensional case in the next section.

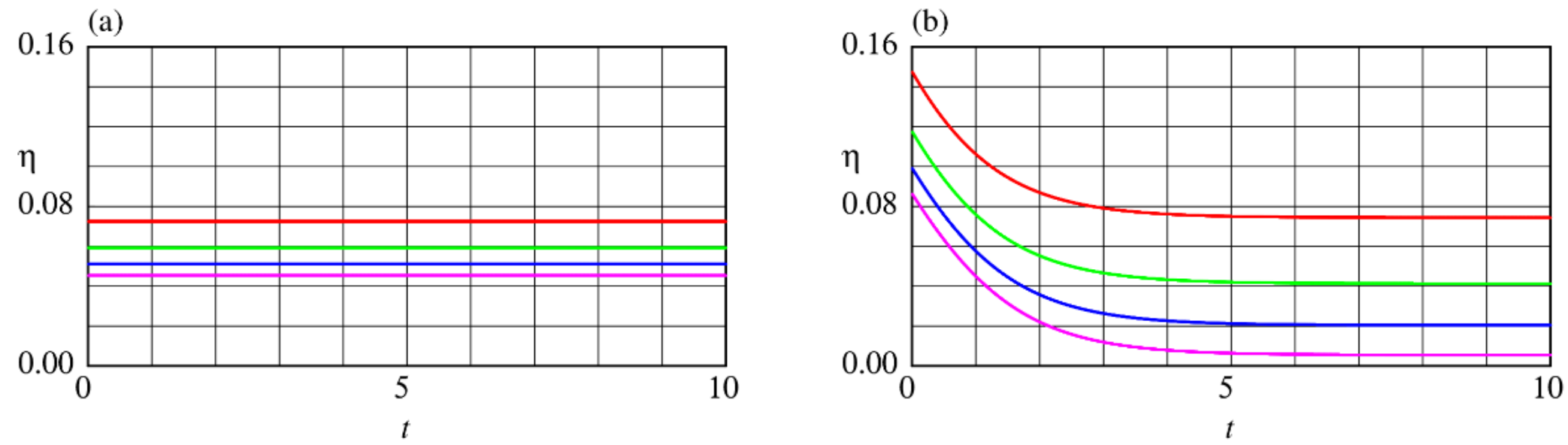

**Figure 1.** Computed time histories of the optimal Lagrangian multiplier $\eta = \eta(\mu_t)$: (a) $U_1$ and (b) $U_2$. The color legends represent the results for $\varepsilon = 0.150$ (red), $\varepsilon = 0.225$ (green), $\varepsilon = 0.300$ (blue), and $\varepsilon = 0.375$ (magenta).

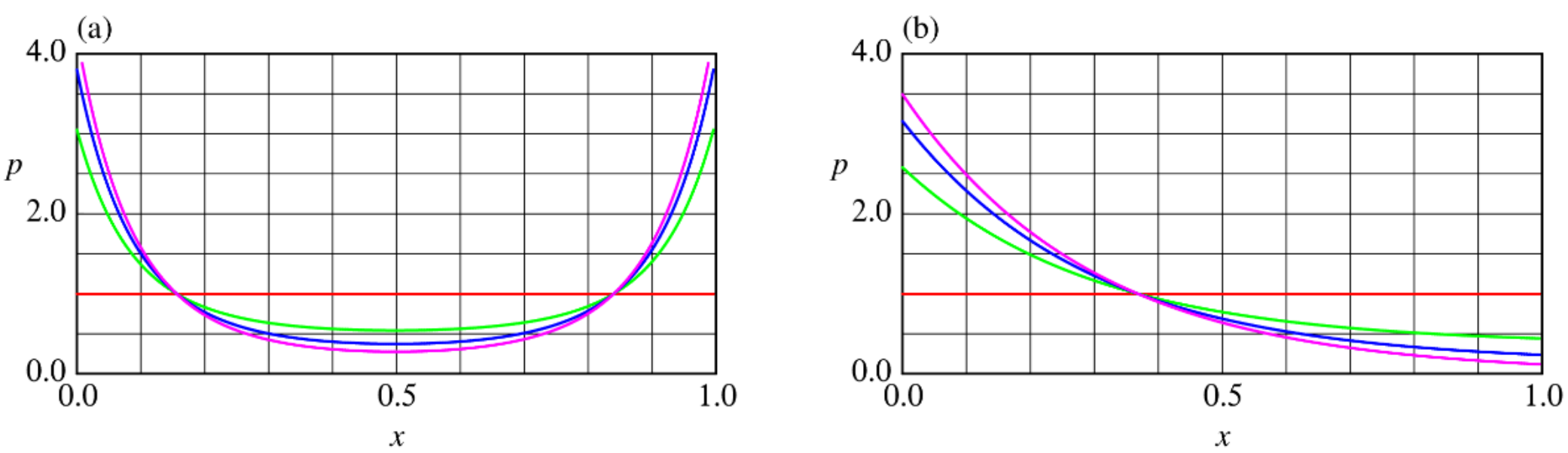

**Figure 2.** Computed time histories of the probability densities $p = p_t(x)$: (a) $U_1$ and (b) $U_2$. The color legends represent the results for $t = 0$ (red), $t = 1$ (green), $t = 2$ (blue), and $t = 10$ (magenta).

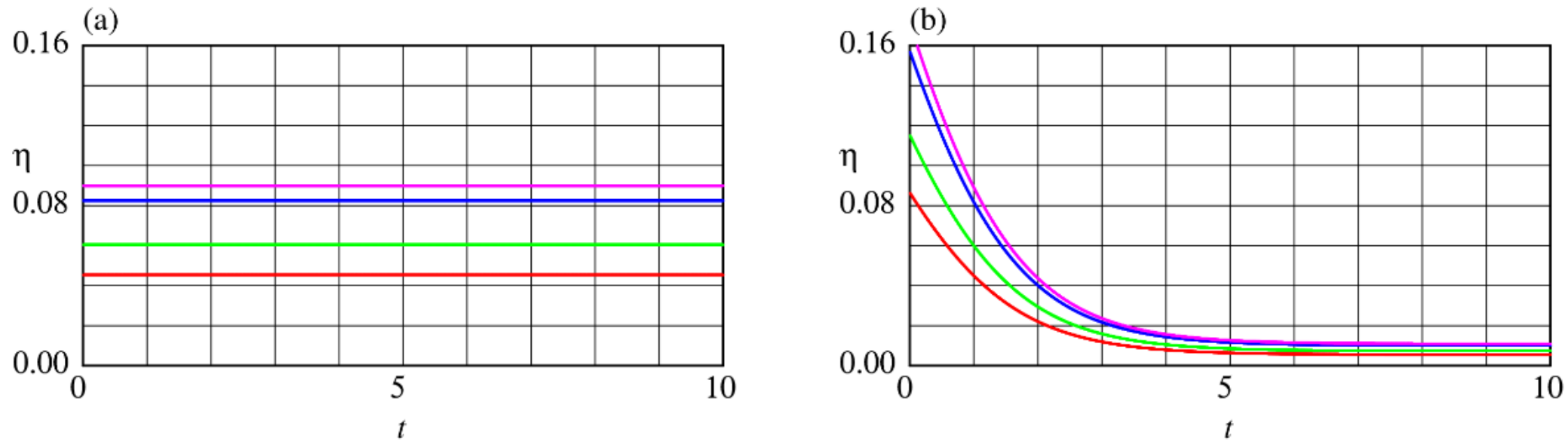

**Figure 3.** Computed time histories of the optimal Lagrangian multiplier $\eta = \eta(\mu_t)$: (a) $U_1$ and (b) $U_2$. The color legends represent the results for $\delta = 1$ (red), $\delta = 2$ (green), $\delta = 10$ (blue), and $\delta = 100$ (magenta).

**Table 1.** Computed average value of agent actions and their errors from the exact value of 0.25. The convergence rate is defined as follows: the convergence rate of the computed $\varepsilon$ at level $I$ against that of $I-1$ is estimated by $\log_{\frac{\varepsilon_I}{\varepsilon_{I-1}}}\left(\frac{\mathrm{Error}_{I-1}}{\mathrm{Error}_I}\right)$.

| $I$ | $\varepsilon$ | Average | Error | Convergence rate |
|---|---|---|---|---|
| 1 | 0.150 | 3.44.E-01 | 9.43.E-02 | |
| 2 | 0.225 | 3.11.E-01 | 6.08.E-02 | 1.08.E+00 |
| 3 | 0.300 | 2.83.E-01 | 3.32.E-02 | 2.10.E+00 |
| 4 | 0.375 | 2.60.E-01 | 9.63.E-03 | 5.56.E+00 |

***Sensitivity against $\chi$***

We also investigate the influences of the regularization in the cost constraint in the BNN and replicator models for the second utility $U = U_2$. Here, we set $\delta = 1$, $\varepsilon = 0.1$, $\chi = 10^{-5}$, and $\xi = 2$ unless otherwise specified. Numerical solutions converge toward equilibrium as time elapses (**Figure 4**; equilibria are discussed in greater detail later); at the same time, the optimal multiplier $\eta$ decreases to some finite value (**Figure 5**). This phenomenon is considered to correspond to a situation where agents are moving toward an equilibrium solution, the search becomes less important, as for the logit model. The same seems not computationally hold true for the unregularized version ($\chi = 0$), which is not covered by **Proposition 1**; the cost $\varepsilon$ is fixed to 0.1, while the optimal multiplier $\eta$ becomes closer to 0 and then oscillates at small values, and eventually the numerical solutions fail. This phenomenon is considered due to $\Phi = U_2$ in theory at a Nash equilibrium; at the equilibrium, we have $\Phi = U_2 = const$, with which the cost constraint becomes $\eta(\mu_t)^{2+\xi} = \chi\varepsilon$ when $\chi > 0$; however, $\eta(\mu_t) = 0$ when $\chi = 0$, with which the quantity $\Phi / \eta$ becomes undetermined and the computation fails.

We further study the regularization method. As discussed above, the regularization term included in the cost constraint plays a role in both mathematical and computational aspects to obtain solutions. Theoretically, $\chi$ is allowed to be an arbitrary positive value, but it should not be too large in computation because an excessively large value of $\chi$ may result in a model whose cost constraint deviates from the original one (24) with $R \equiv 0$. **Figure 6** shows the normalized true exploration cost given by

$$E_t = 1 - \frac{\chi}{\varepsilon \eta(\mu_t)^{2+\xi}} \geq 0, \tag{49}$$

which equals $E_t = 1$ when $\chi = 0$. Replicator models have smaller values of $E_t$ than the BNN model, suggesting faster convergence to equilibria. For both models, $E$ is smaller at earlier times for smaller $\chi$ values corresponding to weaker regularization effects and decays to 0 as time elapses. An unresolved question is whether $\eta(\mu_t) = 0$ occurs in a finite time or not and its dependence on $\chi$, which would need another theoretical study to investigate EGDs with singular coefficients. **Figure 7** compares time histories of the optimal Lagrangian $\eta$ and true exploration cost for different values of the power $\xi$ of regularization. Specifying a larger value of $\xi$ corresponding to stronger regularization results in faster convergence particularly for earlier times; however, truly faster convergence is achieved at smaller $\xi$; for the case $\xi = 0$, the convergence to an equilibrium is computationally achieved at time 5.8.

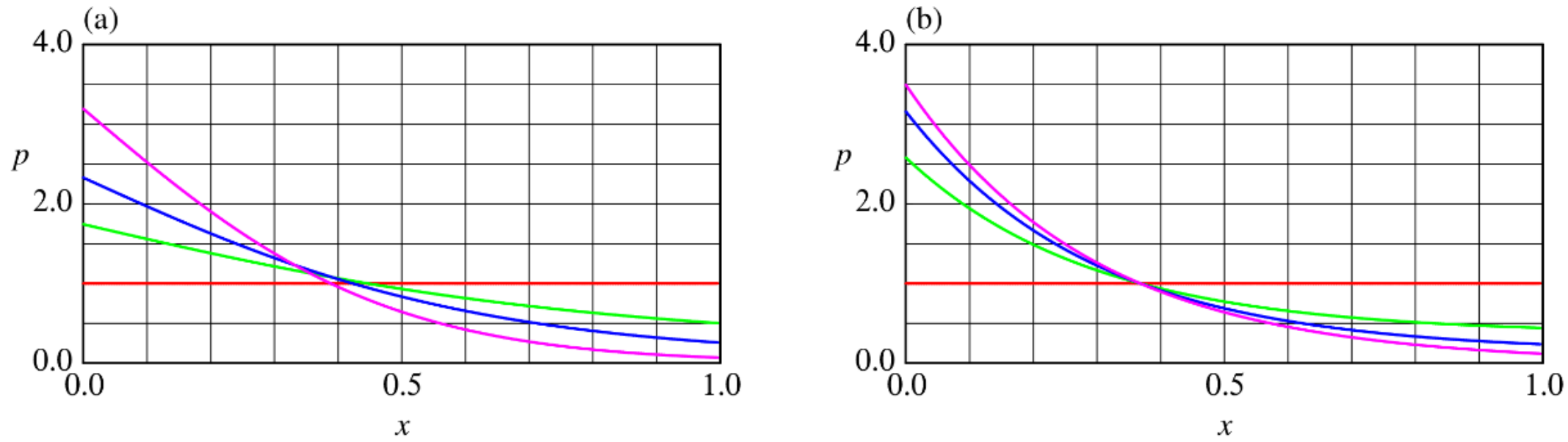

**Figure 4.** Computed time histories of the probability densities $p = p_t(x)$: (a) BNN model and (b) replicator model. The color legends represent the results for $t = 0$ (red), $t = 1$ (green), $t = 2$ (blue), and $t = 10$ (magenta). The computed probability densities are close to equilibria at $t = 10$.

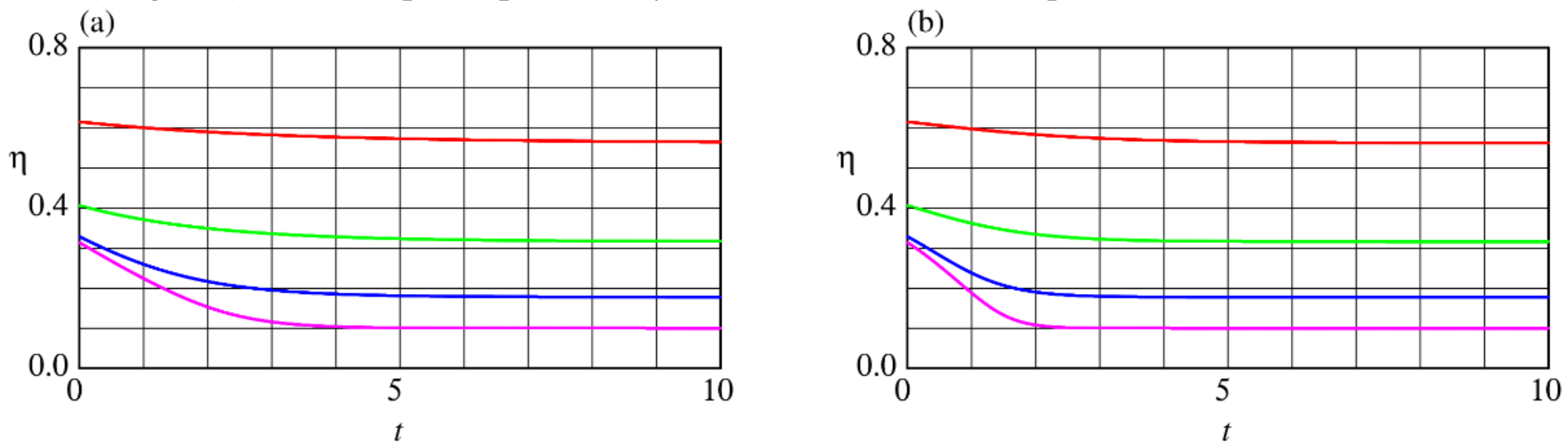

**Figure 5.** Computed time histories of the optimal Lagrangian multiplier $\eta = \eta(\mu_t)$: (a) BNN model and (b) replicator model. The color legends represent the results for $\chi = 10^{-2}$ (red), $\chi = 10^{-3}$ (green), $\chi = 10^{-4}$ (blue), and $\chi = 10^{-5}$ (magenta).

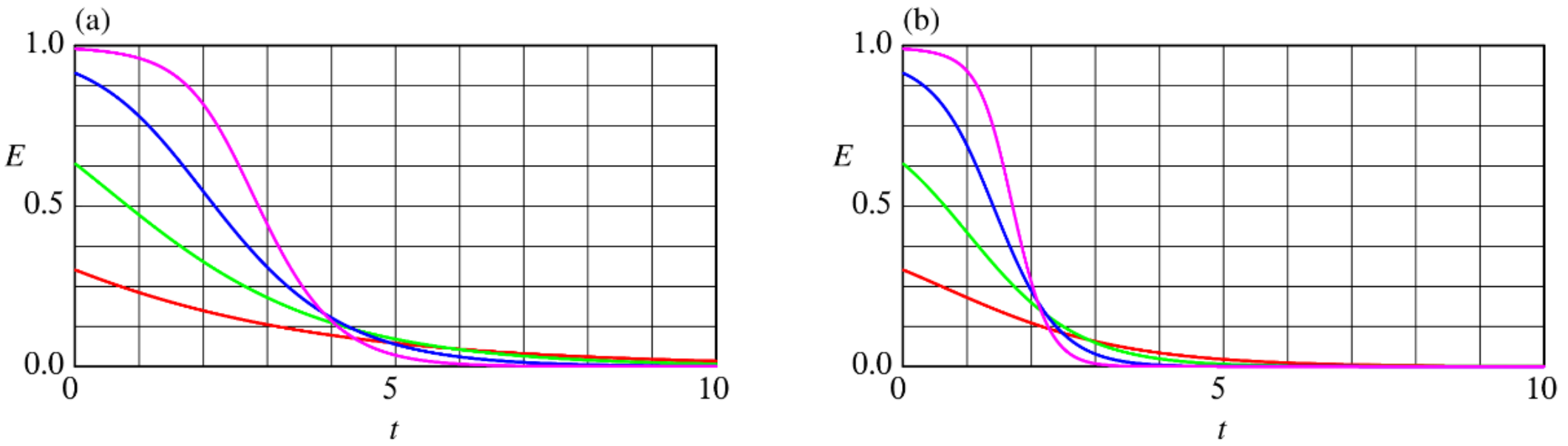

**Figure 6.** Computed time histories of the true exploration cost $E = E_t$: (a) BNN model and (b) replicator model. The color legends represent the results for $\chi = 10^{-2}$ (red), $\chi = 10^{-3}$ (green), $\chi = 10^{-4}$ (blue), and $\chi = 10^{-5}$ (magenta).

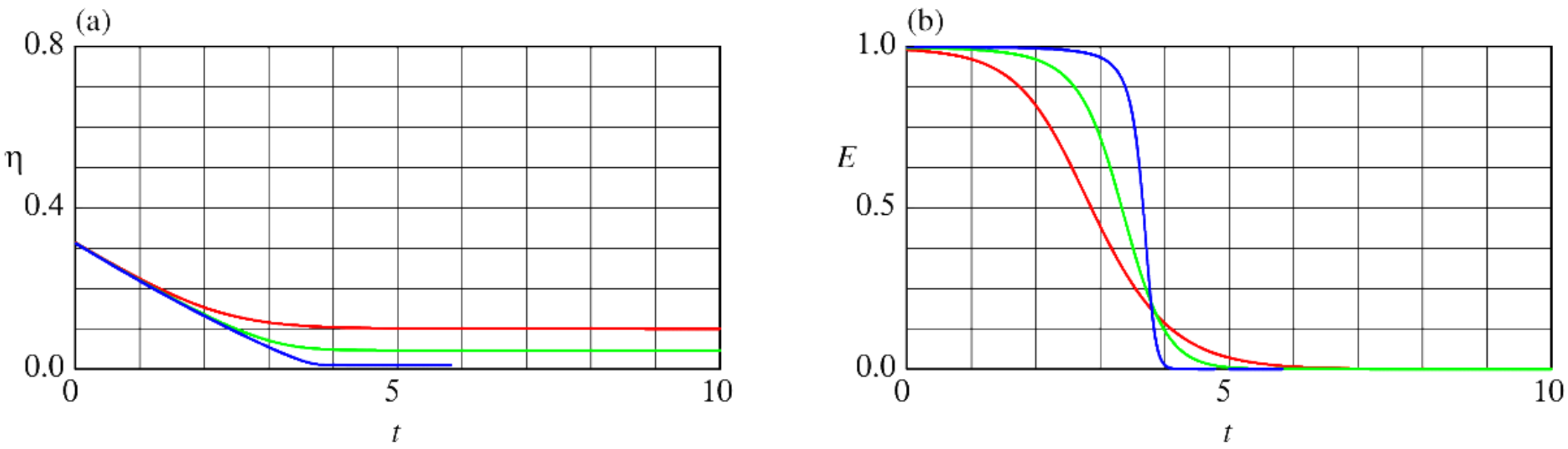

**Figure 7.** Computed time histories of the (a) optimal Lagrangian multiplier $\eta = \eta(\mu_t)$ and (b) true exploration cost $E = E_t$ for the BNN model against different values of $\xi$: $\xi = 2$ (red), $\xi = 1$ (green), $\xi = 0$ (blue).

*Equilibria*

The (Nash) equilibria of the second utility $U = U_2$ are not unique, and here, we compare computed equilibria for each model. We examine the cases $\delta = 1$ and $\delta = 10^8$, where the latter is the myopic case with approximately $U = \Phi$. We fix $\varepsilon = 0.375$ with which numerical solutions by the logit model give an average reasonably close to 0.25 of the theoretical Nash equilibria according to **Table 1**. The computed stationary probability densities for the three models are compared in **Figure 8**, and the results suggest that the equilibria selected by the three models differ from each other, suggesting that different protocols possibly yield distinctive equilibria in our EGDs for both large and small discount rates.

The average of agent actions is very close to 0.25 (with error levels smaller than 0.0001%) for both $\delta = 1$ and $\delta = 10^8$, although the corresponding equilibria differ according to **Figure 8**. Therefore, for the EGDs of the BNN and replicator types with finite $\delta$, the Nash equilibria are attained despite the forward-looking EGD (26) using $\Phi$ instead of $U$, which is considered because both are close to 0, as already reported above. For these two models, a larger $\delta$ results in a probability density with higher maxima at $x = 0$, meaning that the total number of agents with small harvesting rates, which would contribute more to sustainable resource management, increases if they have a longer perspective.

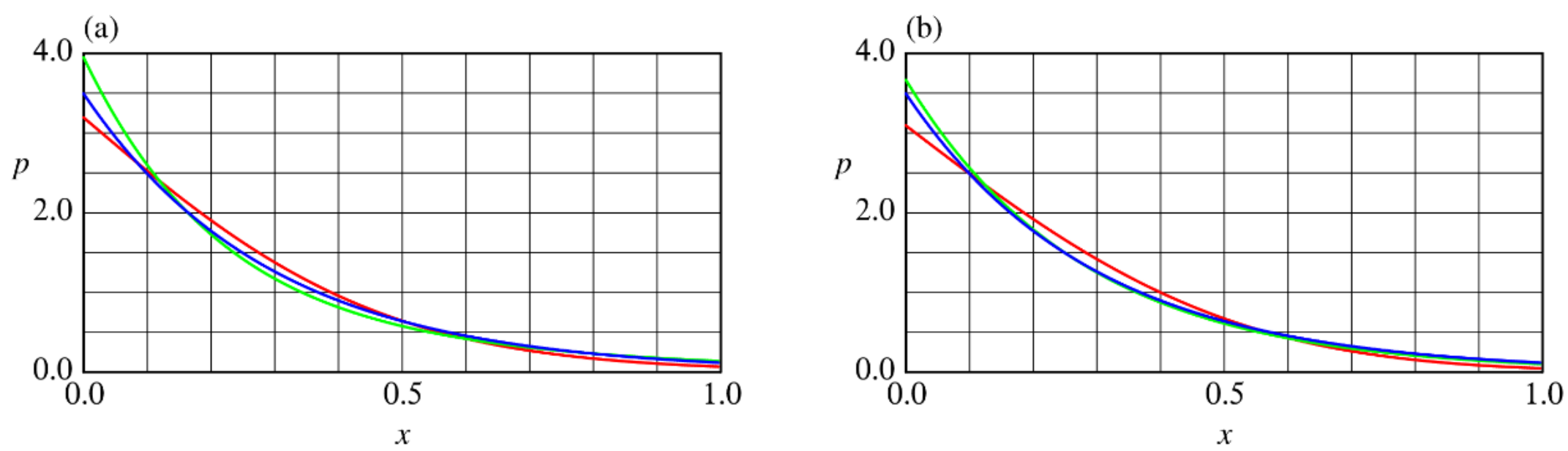

**Figure 8.** Computed stationary probability densities $p = p(x)$ for different models: (a) $\delta = 1$ model and (b) $\delta = 10^8$. The color legends represent the results for the BNN model (red), logit model (green), and replicator model (blue).

### 4.2 Two-dimensional problem

#### 4.2.1 Problem setting

We consider a two-dimensional extension of our EGD focusing on the logit model, where the domain $\Omega$ of agent actions is given by the cube $[0,1]$ in an $x$-z plane and the utility $U_3$, which is based on $U_2$, is set as

$$U_3(x,z,\mu) = \left(h(z) f\left(\int_\Omega y\mu(\mathrm{d}y,\mathrm{d}z)\right) - c\right)x \tag{50}$$

with a function $h:[0,1]\to\mathbb{R}_+$. Now, the action contains the two independent variables $(x,z)$, where $x$ represents the harvesting intensity of the resource and $z$ represents its efficiency. We assume that $h$ is strictly increasing. Under this problem setting, agents can choose not only how much to harvest but also how to do so. These utility models add heterogeneity, which is the new degree-of-freedom $z$, to model different functional shapes of utilities, and this $z$ can be updated as for $x$. In the present case, an efficient harvesting method would yield higher profits for some individuals with high $z$ because of the increasing nature of $h$, while at the same time, it possibly leads to a decrease in utility on average. These mechanisms govern how agents update their two-dimensional actions as time elapses and what equilibrium point they will reach.

Existing studies on common-pool resource management based on EGDs and related models with heterogeneous agents include pairwise interactions on networks [37,44], agents on lattices interacting through heterogeneous utilities [46], ODE models accounting for the resource growth rate [35] and those considering maladaptation under climate change [3], all highlighting the relevance of modeling the social phenomena of heterogeneous agents. To the best of the author's knowledge, the existing studies do not thoroughly address common-pool resource management by using EGDs with mean-field structures. In this view, the computational analysis here is still germinating but provides insights into more advanced modeling for common-pool resource management.

The mathematical analysis results obtained in earlier sections carry over to the two-dimensional case. **Proposition 2** applies to this two-dimensional case if the Lipschitz continuity of utilities, which is **(A4)**, is imposed for both the $x$ direction and the $z$ direction. The EGD is discretized in both the $x$ and $z$ directions by straightforwardly extending the one-dimensional numerical method (**Section A.3**). The discretization parameter values are set as $\Delta t = 0.005$ and $\Delta x = \Delta z = 1/N = 1/250$. We set $h(z) = z$. The values of the other parameters in the utility are the same as those for $U_2$. Unless otherwise specified, we set $\delta = 1$ and $\varepsilon = 1.5$.

#### 4.2.2 Computational results and discussion

The time history of the probability density $p$, which is given by $\mu_{i,j,k}/(\Delta x \Delta z)$ at each time instance and cell, is shown in **Figure 9**. A profile concentrated around the boundaries $x = 0$ and $z = 1$ is developed, and the numerical solution is almost stationary at $t = 10$. Unlike in the one-dimensional case, the probability density of agent actions does not necessarily decrease with increasing $x$ and increases with increasing $z$. Because the initial condition is a uniform distribution on $\Omega$, the computational results suggest the coexistence of agents with (small $z$) and without sustainability concerns (large $z$).

**Figure 10** compares the value functions $\Phi$ for $\delta = 1,\ 10^8$ at stationary states, demonstrating that the two cases have different value functions with each other up to truncation errors. A key to understanding this phenomenon is the formula (108) for the one-dimensional case in **Section A.3** (and the two-dimensional version examined here); this formula, along with its unique solvability of the EGD, yields the following equality: For $\delta = \delta_1,\ \delta_2$, it follows that

$$\frac{\delta_1}{\delta_1 + 1} \frac{U(x,\mu_t)}{\eta(\mu_t)\big|_{\delta=\delta_1}} = \frac{\delta_2}{\delta_2 + 1} \frac{U(x,\mu_t)}{\eta(\mu_t)\big|_{\delta=\delta_2}}, \tag{51}$$

showing that if $U(x,\mu_t) \neq 0$, then $(\delta+1)\eta/\delta$ does not depend on $\delta$. This further implies that the logit function equals for different values of $\delta$ by (106) in **Section A.3**. The computational results obtained here imply that the same observation applies to the discrete case. In this view, the discount rate and cost constraint play essentially the same role, and our numerical solutions correctly capture this property. Note that this property does not apply to the BNN and replicator models, as demonstrated in **Figure 8**.

Finally, the stationary probability densities for different values of $\varepsilon$ are shown in **Figure 11**. The computed probability densities decrease in $x$ for small values of $\varepsilon$, while they are not monotonic for larger $\varepsilon$, as also shown in **Figure 9**. This monotonicity, and in particular, the concentration of probability densities along the boundary $z = 1$, becomes more significant as $\varepsilon$ increases. **Figure 12** shows the computed time histories of the optimal Lagrangian multiplier $\eta$ for different values of $\varepsilon$. As in the one-dimensional case, $\eta$ decreases in $\varepsilon$, suggesting that the two-dimensional case inherits the perspectives of the agents in the one-dimensional case. Moreover, a faster convergence to a nearly constant $\eta$ is achieved for a larger $\varepsilon$, suggesting that paying a larger cost enables agents to find an equilibrium faster.

As demonstrated in this paper, our framework for forward-looking EGDs applies to both one- and multidimensional cases. We conjecture that the same will be true for higher dimensions, but tackling such a problem would require another computational approach, which is currently under investigation by the author(s).

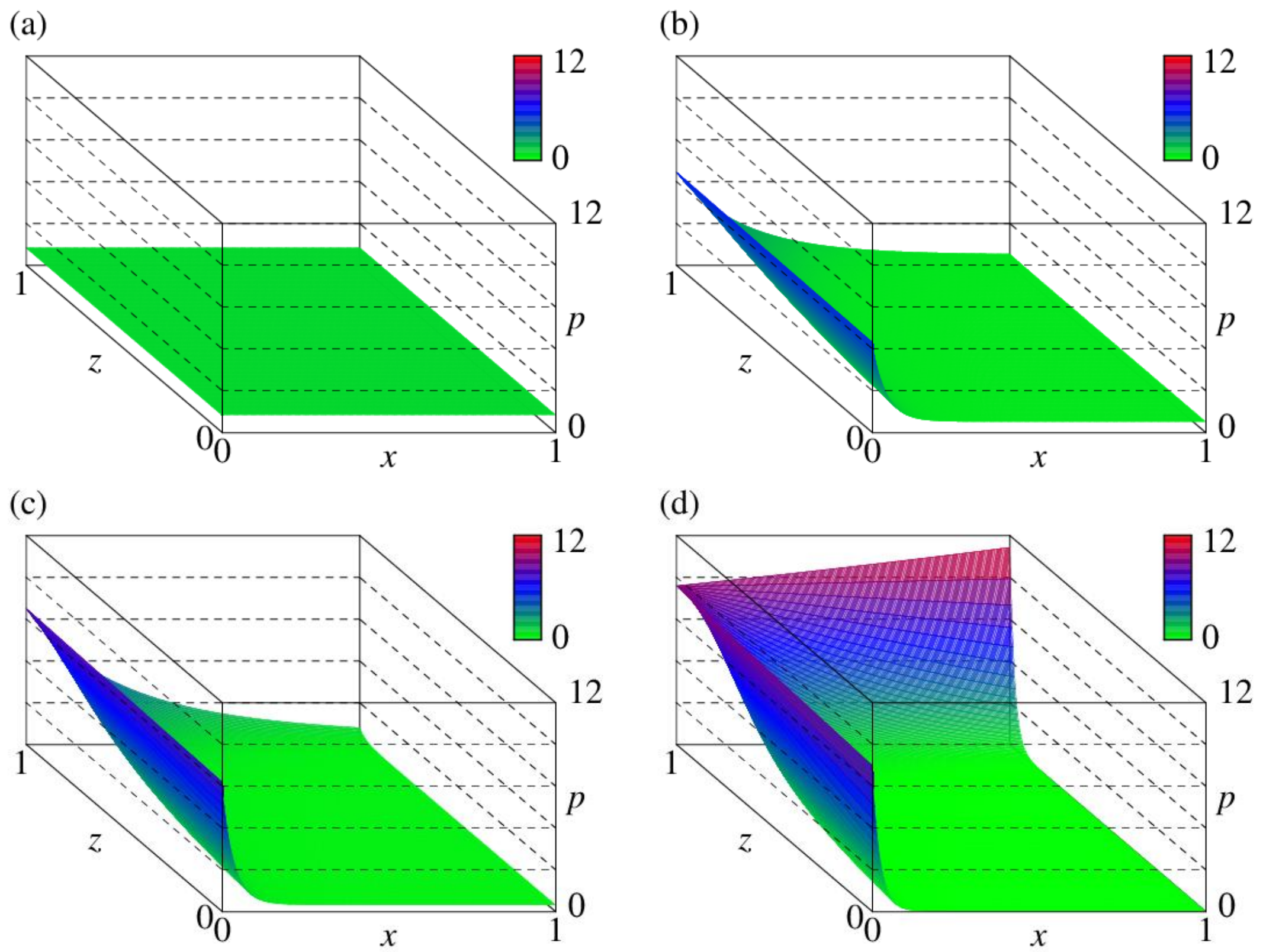

**Figure 9.** Computed time histories of the probability densities $p = p_t(x,z)$: (a) $t = 0$, (b) $t = 0.5$, (c) $t = 1$, and (d) $t = 5$.

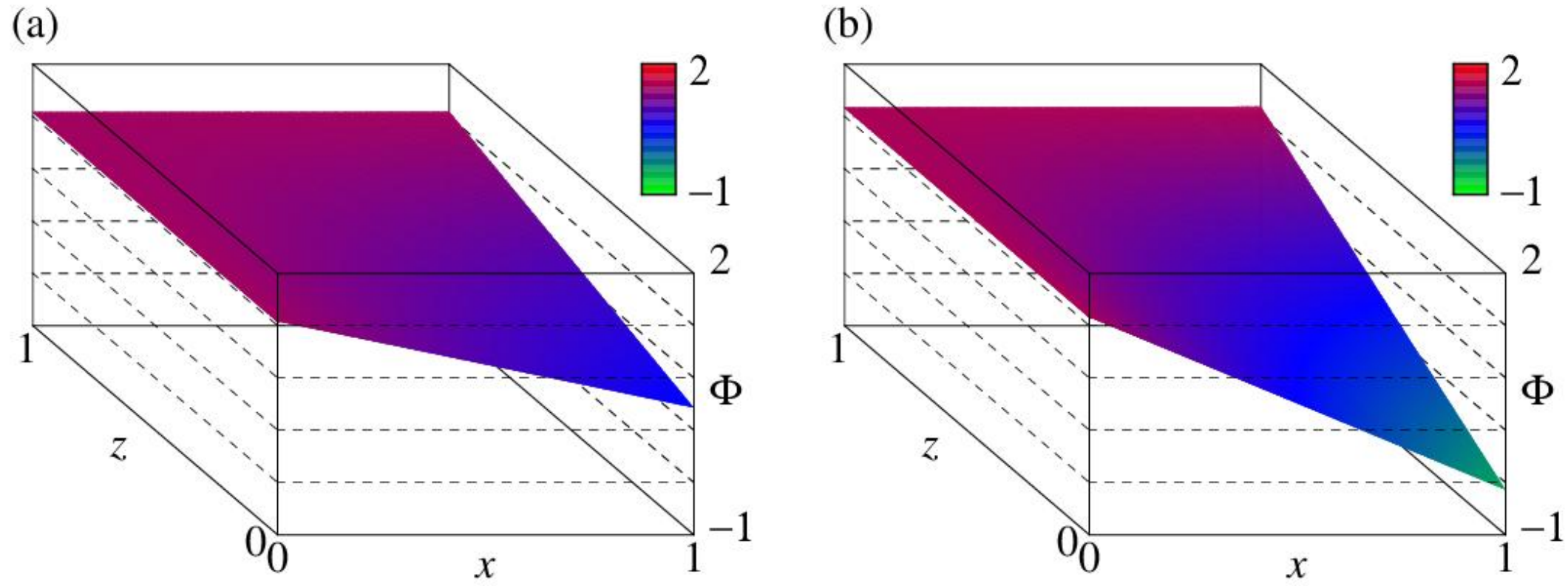

**Figure 10.** Computed value functions $\Phi$ at stationary states: (a) $\delta = 1$ and (b) $\delta = 10^8$.

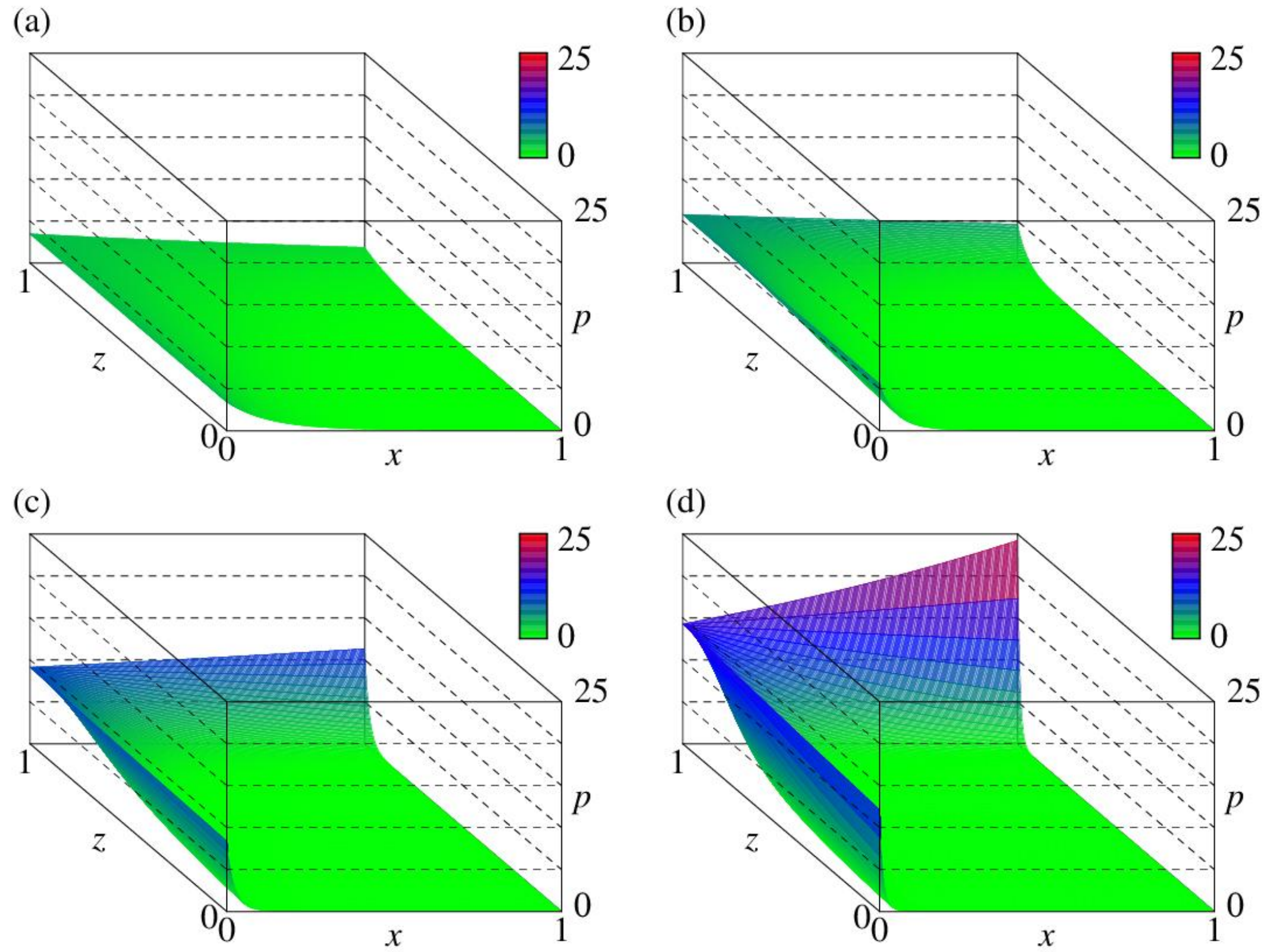

**Figure 11.** Computed stationary probability densities $p$: (a) $\varepsilon = 0.5$, (b) $\varepsilon = 1.0$, (c) $\varepsilon = 1.5$, and (d) $\varepsilon = 2.0$.

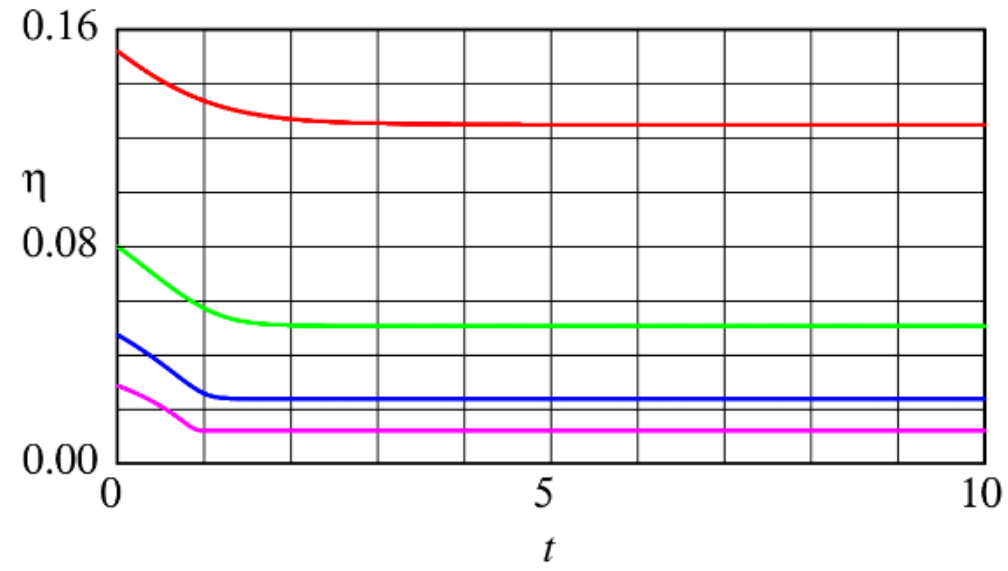

**Figure 12.** Computed time histories of the optimal Lagrangian multiplier $\eta = \eta(\mu_t)$ in the two-dimensional problem. The color legends represent the results for $\varepsilon = 0.5$ (red), $\varepsilon = 1.0$ (green), $\varepsilon = 1.5$ (blue), and $\varepsilon = 2.0$ (magenta).

## 5. Conclusion

We proposed a novel EGD, called forward-looking EGD, in a one-dimensional space that accounts for forward-looking behavior and exploration cost. The resulting system contains a PIDE of the time-dependent probability measure of agent actions and a static HJB equation of the value function. By focusing on the BNN, replicator, and logit protocols, we proved the unique existence for our forward-looking EGD under certain conditions of utility and discount rate, and regularization. We also presented a numerical method for discretizing the EGD and applied it to one-dimensional cases and a two-dimensional problem of common-pool resource management. This paper thus addressed both the theoretical and application aspects of a new mathematical model.

This paper did not address forward-looking EGDs with time-dependent coefficients, which probably need to use a time-dependent HJB equation to which a Banach fixed-point theorem may not apply but Schauder one does [4]. In this case, the existence of solutions can be proven, but uniqueness will not be straightforward. Another remaining issue is the extension of the obtained mathematical results to general pairwise comparison protocols that should be formulated in a more abstract way. In this regard, **Assumption 1**, particularly **(A1) and (A5)**, can be improved depending on the problem at hand.

From an application side, EGDs with higher dimensions will appear in engineering problems, which need to be computed by a data-driven numerical method [52]. Common-pool resource management with spatial heterogeneity [32] would also be a viable next step of our study toward the development of a more realistic EGD. Another potentially interesting topic is to extend the proposed framework to EGDs that were not covered by this study, such as the nonlocal cubic model [8]. These issues are currently under investigation by the author(s) to develop advanced models for sustainable resource management.

## Appendix

### A.1 Nash equilibria

We define Nash equilibria following Section 3 in Cheung and Lahkar [13], which is an equilibrium in which any agent does not gain any advantage even if he/she changes action unless the other agents change their actions.

***Definition A1***

*For a utility* $U:\Omega\times\mathcal{P}\to\mathbb{R}$, $\mu^{*}\in\mathcal{P}$ *with support* $S\left(\mu^{*}\right)$ *is said to be a Nash equilibrium if the following condition is satisfied:*

$$U\left(y,\mu^{*}\right)\le U\left(z,\mu^{*}\right)\ \textit{for all}\ \ z\in S\left(\mu^{*}\right)\ \ \textit{and all}\ \ y\in\Omega. \tag{52}$$

### A.2 Proofs

***Proof of Proposition 1***

First, in view of Theorem 1 in Cheung [12], the main task is to show that the quantity $\eta\left(\sigma\right)^{-1}\Phi\left(x,\sigma\right)$ at each $x\in\Omega$ is bounded and Lipschitz continuous with respect to an arbitrary $\sigma\in\mathcal{M}_2$. Additionally, we need to show the unique existence of the couple $\left(\eta\left(\sigma\right),\Phi\left(\cdot,\sigma\right)\right)\in\mathbb{R}_{+}\times C$ for each $\sigma\in\mathcal{M}_2$ and their uniform bounds. For each $\tau\in\mathcal{M}_2$, its positive part is denoted as $\tau^{+}\in\mathcal{M}_2$. We have $\left\|\tau^{+}\right\|_{\mathrm{TV}}\le\left\|\tau\right\|_{\mathrm{TV}}\le 2$. For each $\sigma,\omega\in\mathcal{M}_2$ and $w\in\left[0,1\right]$, we set $\theta=\left(1-w\right)\kappa+w\sigma\in\mathcal{M}_2$. The proof is divided into several parts.

***Step 1. Auxiliary HJB equation: unique existence***

We fix arbitrary $m>0$. We fix arbitrary $\sigma,\omega\in\mathcal{M}_2$ and $w\in\left[0,1\right]$ and set $\theta=w\kappa+\left(1-w\right)\sigma\in\mathcal{M}_2$. We show that the following auxiliary HJB equation yields a unique solution $\Psi\in C$ such that $0\le\Psi\le\bar{U}$ when $\delta$ is large:

$$\delta\Psi\left(x\right)=\delta U\left(x,\sigma\right)+\sup_{u_x\left(\cdot\right)\ge 0}\left\{\int_{\Omega}u_x\left(z\right)\left(\Psi\left(z\right)-\Psi\left(x\right)\right)\theta^{+}\left(\mathrm{d}z\right)-m\int_{\Omega}\frac{1}{2}\left(u_x\left(z\right)\right)^{2}\theta^{+}\left(\mathrm{d}z\right)\right\},\ \ x\in\Omega. \tag{53}$$

First, if $\Psi\in C$, then the second term on the right-hand side of (53) is rewritten as follows:

$$\sup_{u_x\left(\cdot\right)\ge 0}\left\{\int_{\Omega}u_x\left(z\right)\left(\Psi\left(z\right)-\Psi\left(x\right)\right)\theta^{+}\left(\mathrm{d}z\right)-m\int_{\Omega}\frac{1}{2}\left(u_x\left(z\right)\right)^{2}\theta^{+}\left(\mathrm{d}z\right)\right\}=\frac{1}{2m}\int_{\Omega}\left(\Psi\left(z\right)-\Psi\left(x\right)\right)_{+}^{2}\theta^{+}\left(\mathrm{d}z\right) \tag{54}$$

since

$$\operatorname*{arg\,max}_{u_x\left(\cdot\right)\ge 0}\left\{\int_{\Omega}u_x\left(z\right)\left(\Psi\left(z\right)-\Psi\left(x\right)\right)\theta^{+}\left(\mathrm{d}z\right)-m\int_{\Omega}\frac{1}{2}\left(u_x\left(z\right)\right)^{2}\theta^{+}\left(\mathrm{d}z\right)\right\}=\left(\frac{\Psi\left(\cdot\right)-\Psi\left(x\right)}{m}\right)_{+}\ \left(:=\hat{u}\left(z\right)\right). \tag{55}$$

Then, the auxiliary HJB equation (53) is rewritten as follows:

$$\Psi(x) = U(x,\sigma) + \frac{1}{2m\delta}\int_{\Omega}\left(\Psi(z) - \Psi(x)\right)_{+}^{2}\theta^{+}(\mathrm{d}z) \quad \left(:= \mathbb{A}(\Psi)(x)\right), \quad x \in \Omega . \tag{56}$$

Any solution $\Psi \in C$ to (56) must satisfy $0 \le \Psi \le \bar{U}$ in $\Omega$. Indeed, if $\Psi$ is minimized at some $\hat{x} \in \Omega$, then

$$\Psi(\hat{x}) \ge 0 + \frac{1}{2m\delta}\int_{\Omega}\left(\Psi(z) - \Psi(\hat{x})\right)_{+}^{2}\theta^{+}(\mathrm{d}z) \ge 0 , \tag{57}$$

yielding $\Psi(\hat{x}) \ge 0$. Similarly, if $\Psi$ is maximized at some $\hat{x} \in \Omega$, then

$$\Psi(\hat{x}) \le \bar{U} + \frac{1}{2m\delta}\int_{\Omega}\left(\Psi(z) - \Psi(\hat{x})\right)_{+}^{2}\theta^{+}(\mathrm{d}z) = \bar{U} + 0 = \bar{U} , \tag{58}$$

yielding $\Psi(\hat{x}) \le \bar{U}$.

We show that $\mathbb{A}: C \to C$ is bounded and strictly contractive in $C$ for a sufficiently large $\delta > 0$. Note that $\mathbb{A}(\Psi) \in C$ when $\Psi \in C$ because of **(A4) in Assumption 1**, i.e., (19). The *a priori* bound derived above implies that it suffices to consider those with $0 \le \Psi \le \bar{U}$ in $\Omega$ (more specifically, we first consider another auxiliary HJB equation with $\Psi$ replaced by $0 \vee \left(\Psi \wedge \bar{U}\right)$ on the right-hand side and prove the unique existence of its solution; after that, the same solution is solved (56); see Proof of Proposition 2 in Yoshioka et al. [49]). In this view, the boundedness is proven as follows: For any $\Psi_1 \in C$ with $0 \le \Psi_1 \le \bar{U}$,

$$\begin{aligned}
\left\|\mathbb{A}(\Psi_1)\right\|_{\infty} &\le \bar{U} + \frac{1}{2m\delta}\int_{\Omega}\left(\bar{U} - 0\right)_{+}^{2}\theta^{+}(\mathrm{d}z) \\
&\le \bar{U} + \frac{\bar{U}^2}{2m\delta}\left\|\theta\right\|_{\mathrm{TV}} \\
&= \bar{U} + \frac{\bar{U}^2}{m\delta}.
\end{aligned} \tag{59}$$

The strict contraction property is proven as follows: For any $\Psi_1, \Psi_2 \in C$ with $0 \le \Psi_1, \Psi_2 \le \bar{U}$, if $\delta > 4m^{-1}\bar{U}$, then

$$
\begin{aligned}
&\left\|\mathbb{A}(\Psi_1)-\mathbb{A}(\Psi_2)\right\|_\infty \\
&=\frac{1}{2m\delta}\left|\int_\Omega\left(\Psi_1(z)-\Psi_1(\cdot)\right)_+^2\theta^+(\mathrm{d}z)\ -\int_\Omega\left(\Psi_2(z)-\Psi_2(\cdot)\right)_+^2\theta^+(\mathrm{d}z)\ \right| \\
&\leq\frac{1}{2m\delta}\int_\Omega\left|\left(\Psi_1(z)-\Psi_1(\cdot)\right)_+^2-\left(\Psi_2(z)-\Psi_2(\cdot)\right)_+^2\right|\theta^+(\mathrm{d}z) \\
&=\frac{1}{2m\delta}\int_\Omega\left|\left(\Psi_1(z)-\Psi_1(\cdot)\right)_++\left(\Psi_2(z)-\Psi_2(\cdot)\right)_+\right|\left|\left(\Psi_1(z)-\Psi_1(\cdot)\right)_+-\left(\Psi_2(z)-\Psi_2(\cdot)\right)_+\right|\theta^+(\mathrm{d}z) \\
&\leq\frac{\bar{U}}{m\delta}\int_\Omega\left|\left(\Psi_1(z)-\Psi_1(\cdot)\right)_+-\left(\Psi_2(z)-\Psi_2(\cdot)\right)_+\right|\theta^+(\mathrm{d}z) \qquad (60)\\
&\leq\frac{\bar{U}}{m\delta}\int_\Omega\left|\left(\Psi_1(z)-\Psi_1(\cdot)\right)-\left(\Psi_2(z)-\Psi_2(\cdot)\right)\right|\theta^+(\mathrm{d}z) \\
&\leq\frac{2\bar{U}}{m\delta}\left\|\Psi_1-\Psi_2\right\|_\infty\left\|\theta\right\|_{\mathrm{TV}} \\
&\leq\frac{4\bar{U}}{m\delta}\left\|\Psi_1-\Psi_2\right\|_\infty \\
&<\left\|\Psi_1-\Psi_2\right\|_\infty .
\end{aligned}
$$

Therefore, by Banach's fixed-point theorem (e.g., Theorem 5.7 in Brezis [9]), if $\delta>4m^{-1}\bar{U}$, then the auxiliary HJB equation admits the unique solution denoted by $\Psi_m\in C$ such that $0\leq\Psi_m\leq\bar{U}$.

***Step 2. Auxiliary HJB equation: cost constraint***

Under the setting of **Step 2**, from (29) and (55), with $\Psi=\Psi_m$ we have

$$
\begin{aligned}
&\int_{x\in\Omega}\int_{z\in\Omega}\frac{1}{2}\left(u_x^*(z)\right)^2\theta^+(\mathrm{d}z)\theta^+(\mathrm{d}x)+\frac{\chi}{m^2} \\
&=\frac{1}{2m^2}\int_{x\in\Omega}\int_{z\in\Omega}\left(\Psi_m(z)-\Psi_m(x)\right)_+^2\theta^+(\mathrm{d}z)\theta^+(\mathrm{d}x)+\frac{\chi}{m^2} \qquad (61)\\
&=\frac{\delta}{m}\int_{x\in\Omega}\left(\Psi_m(x)-U(x,\sigma)\right)\theta^+(\mathrm{d}x)+\frac{\chi}{m^2} \\
&\left(:=g(m)\right),
\end{aligned}
$$

which is due to the following equality obtained from (56):

$$
\int_{x\in\Omega}\left(\Psi_m(x)-U(x,\sigma)\right)\theta^+(\mathrm{d}x)=\frac{m}{\delta}\frac{1}{2m^2}\int_{x\in\Omega}\int_{z\in\Omega}\left(\Psi_m(z)-\Psi_m(x)\right)_+^2\theta^+(\mathrm{d}z)\theta^+(\mathrm{d}x). \qquad (62)
$$

We show that exactly one $m=m(\sigma)>0$ exists such that

$$
g(m)=\varepsilon . \qquad (63)
$$

Note that $g(m)>0$ ($m>0$) because $\chi>0$. Differentiating $g(m)$ with respect to $m>0$ yields

$$
\begin{aligned}
\frac{\mathrm{d}g(m)}{\mathrm{d}m}&=-\frac{\delta}{m^2}\int_{x\in\Omega}\left(\Psi_m(x)-U(x,\sigma)\right)\theta^+(\mathrm{d}x)+\frac{\delta}{m}\int_{x\in\Omega}\frac{\partial\Psi_m(x)}{\partial m}\theta^+(\mathrm{d}x)-\frac{2\chi}{m^3} \\
&=-\frac{1}{m}\left(\varepsilon-\frac{\chi}{m^2}\right)+\frac{\delta}{m}\int_{x\in\Omega}\frac{\partial\Psi_m(x)}{\partial m}\theta^+(\mathrm{d}x)-\frac{2\chi}{m^3} \qquad (64)\\
&=-\frac{1}{m}\varepsilon-\frac{\chi}{m^3}+\frac{\delta}{m}\int_{x\in\Omega}\frac{\partial\Psi_m(x)}{\partial m}\theta^+(\mathrm{d}x).
\end{aligned}
$$

Here, from (56), we obtain

$$
\begin{aligned}
\frac{\partial \Psi_m(x)}{\partial m} &= \frac{1}{2m\delta}\int_\Omega \frac{\partial}{\partial m}\left(\Psi_m(z)-\Psi_m(x)\right)_+^2 \theta^+(\mathrm{d}z) - \frac{1}{2m^2\delta}\int_\Omega \left(\Psi_m(z)-\Psi_m(x)\right)_+^2 \theta^+(\mathrm{d}z) \\
&= \frac{1}{m\delta}\int_\Omega \left(\frac{\partial \Psi_m(z)}{\partial m}-\frac{\partial \Psi_m(x)}{\partial m}\right)\left(\Psi_m(z)-\Psi_m(x)\right)_+ \theta^+(\mathrm{d}z) \\
&\quad -\frac{1}{2m^2\delta}\int_\Omega \left(\Psi_m(z)-\Psi_m(x)\right)_+^2 \theta^+(\mathrm{d}z)
\end{aligned}
\quad , \ x\in\Omega, \tag{65}
$$

where $\Psi_m \in C$ is considered given here. For (56), by an argument analogous to that in **Step 1**, the maximum of $\frac{\partial \Psi_m}{\partial m}$ in $\Omega$ is not larger than that on the right-hand side of the following inequality:

$$
-\frac{1}{2m^2\delta}\int_\Omega \left(\Psi_m(z)-\Psi_m(x)\right)_+^2 \theta^+(\mathrm{d}z) \le 0, \quad x\in\Omega . \tag{66}
$$

Similarly, the minimum $\frac{\partial \Psi_m}{\partial m}$ on $\Omega$ is not smaller than that on the right-hand side of the following inequality:

$$
-\frac{1}{2m^2\delta}\int_\Omega \left(\Psi_m(z)-\Psi_m(x)\right)_+^2 \theta^+(\mathrm{d}z) \ge -\frac{\bar{U}^2}{2m^2\delta}, \quad x\in\Omega . \tag{67}
$$

We therefore obtain

$$
-\frac{\bar{U}^2}{2m^2\delta} \le \frac{\partial \Psi_m(x)}{\partial m} \le 0, \quad x\in\Omega \tag{68}
$$

if $\frac{\partial \Psi_m}{\partial m}$ exists and belongs to $C$. Since $\Psi_m \in C$ is considered given here, then (65) is considered a linear equation of $\frac{\partial \Psi_m}{\partial m}$ whose unique existence in $C$ follows if $\delta > 4m^{-1}\bar{U}$; we can use an argument analogous to that in **Step 1** for the following inequality: For any $\Xi_1, \Xi_2 \in C$,

$$
\begin{aligned}
&\left|
\begin{aligned}
&\frac{1}{m\delta}\int_\Omega \left(\Xi_1(z)-\Xi_1(z)\right)\left(\Psi_m(z)-\Psi_m(x)\right)_+ \theta^+(\mathrm{d}z) \\
&-\frac{1}{m\delta}\int_\Omega \left(\Xi_2(z)-\Xi_2(z)\right)\left(\Psi_m(z)-\Psi_m(x)\right)_+ \theta^+(\mathrm{d}z)
\end{aligned}
\right| \\
&\le \frac{\bar{U}}{m\delta}\int_\Omega \left|\left(\Xi_1(z)-\Xi_1(z)\right)-\left(\Xi_2(z)-\Xi_2(z)\right)\right| \theta^+(\mathrm{d}z) \\
&\le \frac{4\bar{U}}{m\delta}\left\|\Xi_1-\Xi_2\right\|_\infty .
\end{aligned}
\tag{69}
$$

Consequently, (64) implies

$$
\frac{\mathrm{d}g(m)}{\mathrm{d}m} = -\frac{1}{m}\varepsilon - \frac{\chi}{m^3} < 0, \tag{70}
$$

and hence $g$ is strictly decreasing for $m>0$, and moreover, we have $\lim_{m\to+0} g(m) = +\infty$ (the function $m^{-1}$ is not integrable near $m=0$) and $\lim_{m\to+\infty} g(m) = 0$. This strict monotonicity of $g$ along with the classical intermediate value theorem shows that exactly one $m = m(\sigma) > 0$ satisfies (63).

Finally, we bound $m(\sigma)$ irrespective of $\sigma$ (and $\theta$). For the upper bound, by the elementary equality $(a-b)_+^2+(b-a)_+^2=(a-b)^2$ ($a,b\in\mathbb{R}$), we have

$$\begin{aligned}
\varepsilon &= \frac{1}{2m(\sigma)^2}\int_{x\in\Omega}\int_{z\in\Omega}\left(\Psi_{m(\sigma)}(z)-\Psi_{m(\sigma)}(x)\right)_+^2\theta^+(\mathrm{d}z)\theta^+(\mathrm{d}x)+\frac{\chi}{m(\sigma)^2} \\
&= \frac{1}{4m(\sigma)^2}\int_{x\in\Omega}\int_{z\in\Omega}\left(\Psi_{m(\sigma)}(z)-\Psi_{m(\sigma)}(x)\right)^2\theta^+(\mathrm{d}z)\theta^+(\mathrm{d}x)+\frac{\chi}{m(\sigma)^2} \\
&\leq \frac{\bar{U}^2}{4m(\sigma)^2}\int_{x\in\Omega}\int_{z\in\Omega}\theta^+(\mathrm{d}z)\theta^+(\mathrm{d}x)+\frac{\chi}{m(\sigma)^2} \\
&\leq \frac{\bar{U}^2}{4m(\sigma)^2}\left(\|\theta\|_{\mathrm{TV}}\right)^2+\frac{\chi}{m(\sigma)^2} \\
&\leq \frac{\bar{U}^2+\chi}{m(\sigma)^2}
\end{aligned} \tag{71}$$

and hence there exists a constant $\bar{\eta}>0$ such that

$$m(\sigma)\leq\left(\frac{\bar{U}^2+\chi}{\varepsilon}\right)^{\frac{1}{2}}\ (:=\bar{\eta}). \tag{72}$$

For the lower bound, we have

$$\varepsilon=\frac{1}{2m(\sigma)^2}\int_{x\in\Omega}\int_{z\in\Omega}\left(\Psi_{m(\sigma)}(z)-\Psi_{m(\sigma)}(x)\right)_+^2\theta^+(\mathrm{d}z)\theta^+(\mathrm{d}x)+\frac{\chi}{m(\sigma)^2}\geq\frac{\chi}{m(\sigma)^2} \tag{73}$$

and hence

$$m(\sigma)\geq\left(\frac{\chi}{\varepsilon}\right)^{\frac{1}{2}}\ (:=\underline{\eta}). \tag{74}$$

Consequently, we obtain the strict bound irrespective of $\sigma$ (and $\delta$) as follows:

$$\underline{\eta}\leq m(\sigma)\leq\bar{\eta}. \tag{75}$$

***Step 3. Auxiliary HJB equation: Lipschitz continuity 1***

We fix arbitrary $m,n>0$. We fix arbitrary $\sigma,\omega\in\mathcal{M}_2$ and $w\in[0,1]$ and set $\theta=w\kappa+(1-w)\sigma\in\mathcal{M}_2$ and $\vartheta=w\kappa+(1-w)\omega\in\mathcal{M}_2$. We write the unique solution to (56) as $\Psi_m(\cdot,\sigma)$. We show its Lipschitz continuity with respect to the second argument. We have

$$
\begin{aligned}
&\left\|\Psi_m(\cdot,\sigma)-\Psi_n(\cdot,\omega)\right\|_\infty \\
&\le \left\|U(\cdot,\sigma)-U(\cdot,\omega)\right\|_\infty+\frac{1}{\delta}\left|\frac{1}{n}-\frac{1}{m}\right|\left|\int_\Omega\left(\Psi_n(z,\sigma)-\Psi_n(\cdot,\sigma)\right)_+^2\vartheta^+(\mathrm{d}z)\right| \\
&\quad+\frac{1}{2\delta m}\left|\int_\Omega\left(\Psi_m(z,\sigma)-\Psi_m(\cdot,\sigma)\right)_+^2\theta^+(\mathrm{d}z)-\int_\Omega\left(\Psi_n(z,\omega)-\Psi_n(\cdot,\omega)\right)_+^2\vartheta^+(\mathrm{d}z)\right| \\
&\le L_U\left\|\sigma-\omega\right\|_{\mathrm{TV}}+\frac{\bar{U}^2}{\delta}\frac{|m-n|}{mn}\left\|\vartheta\right\|_{\mathrm{TV}} \\
&\quad+\frac{1}{2\delta m}\left|\int_\Omega\left(\Psi_m(z,\sigma)-\Psi_m(\cdot,\sigma)\right)_+^2\theta^+(\mathrm{d}z)-\int_\Omega\left(\Psi_n(z,\omega)-\Psi_n(\cdot,\omega)\right)_+^2\vartheta^+(\mathrm{d}z)\right| \\
&\le L_U\left\|\sigma-\omega\right\|_{\mathrm{TV}}+\frac{2\bar{U}^2}{\delta}\frac{|m-n|}{mn} \\
&\quad+\frac{1}{2\delta m}\left(\begin{aligned}&\left|\int_\Omega\left(\Psi_m(z,\sigma)-\Psi_m(\cdot,\sigma)\right)_+^2\theta^+(\mathrm{d}z)-\int_\Omega\left(\Psi_n(z,\omega)-\Psi_n(\cdot,\omega)\right)_+^2\theta^+(\mathrm{d}z)\right| \\ &+\left|\int_\Omega\left(\Psi_n(z,\omega)-\Psi_n(\cdot,\omega)\right)_+^2\theta^+(\mathrm{d}z)-\int_\Omega\left(\Psi_n(z,\omega)-\Psi_n(\cdot,\omega)\right)_+^2\vartheta^+(\mathrm{d}z)\right|\end{aligned}\right).
\end{aligned}
\tag{76}
$$

The first integral in the last line of (76) is evaluated as follows:

$$
\begin{aligned}
&\left|\int_\Omega\left(\Psi_m(z,\sigma)-\Psi_m(\cdot,\sigma)\right)_+^2\theta^+(\mathrm{d}z)-\int_\Omega\left(\Psi_n(z,\omega)-\Psi_n(\cdot,\omega)\right)_+^2\theta^+(\mathrm{d}z)\right| \\
&\le\int_\Omega\left|\left(\Psi_m(z,\sigma)-\Psi_m(\cdot,\sigma)\right)_+^2-\left(\Psi_n(z,\omega)-\Psi_n(\cdot,\omega)\right)_+^2\right|\theta^+(\mathrm{d}z) \\
&\le 2\bar{U}\int_\Omega\left|\left(\Psi_m(z,\sigma)-\Psi_m(\cdot,\sigma)\right)_+-\left(\Psi_n(z,\omega)-\Psi_n(\cdot,\omega)\right)_+\right|\theta^+(\mathrm{d}z) \\
&\le 4\bar{U}\left\|\theta\right\|_{\mathrm{TV}}\left\|\Psi_m(\cdot,\sigma)-\Psi_n(\cdot,\omega)\right\|_\infty \\
&\le 8\bar{U}\left\|\Psi_m(\cdot,\sigma)-\Psi_n(\cdot,\omega)\right\|_\infty .
\end{aligned}
\tag{77}
$$

Similarly, the second integral in the last line of (76) is evaluated as follows:

$$
\begin{aligned}
&\left|\int_\Omega\left(\Psi_n(z,\omega)-\Psi_n(\cdot,\omega)\right)_+^2\theta^+(\mathrm{d}z)-\int_\Omega\left(\Psi_n(z,\omega)-\Psi_n(\cdot,\omega)\right)_+^2\vartheta^+(\mathrm{d}z)\right| \\
&\le\left|\int_\Omega\left(\Psi_n(z,\omega)-\Psi_n(\cdot,\omega)\right)_+^2\left(\theta^+(\mathrm{d}z)-\vartheta^+(\mathrm{d}z)\right)\right| \\
&\le\bar{U}^2\left\|\theta-\vartheta\right\|_{\mathrm{TV}} \\
&\le\bar{U}^2\left\|\sigma-\omega\right\|_{\mathrm{TV}} ,
\end{aligned}
\tag{78}
$$

where $\left\|\theta-\vartheta\right\|_{\mathrm{TV}}\le\left\|\sigma-\omega\right\|_{\mathrm{TV}}$. Consequently, if $\delta>4m^{-1}\bar{U}$, then

$$
\left\|\Psi_m(\cdot,\sigma)-\Psi_n(\cdot,\omega)\right\|_\infty\le L_U\left\|\sigma-\omega\right\|_{\mathrm{TV}}+\frac{2\bar{U}^2}{\delta}\frac{|m-n|}{mn}+\frac{4\bar{U}}{\delta m}\left\|\Psi_m(\cdot,\sigma)-\Psi_n(\cdot,\omega)\right\|_\infty+\frac{\bar{U}^2}{2\delta m}\left\|\theta-\vartheta\right\|_{\mathrm{TV}} ,
\tag{79}
$$

and hence

$$
\left\|\Psi_m(\cdot,\sigma)-\Psi_n(\cdot,\omega)\right\|_\infty\le\left(1-\frac{4\bar{U}}{\delta m}\right)^{-1}\left\{\left(L_U+\frac{\bar{U}^2}{2\delta m}\right)\left\|\sigma-\omega\right\|_{\mathrm{TV}}+\frac{2\bar{U}^2}{\delta}\frac{|m-n|}{mn}\right\}.
\tag{80}
$$

***Step 4. Optimal Lagrangian multiplier: Lipschitz continuity***

Under the setting of **Step 3**, we explicitly write the dependence of $g$ on $\sigma$ (resp., $\omega$) as $g_\sigma$ (resp., $g_\omega$). By (63), we have

$$
m(\sigma)^2=\frac{1}{2\varepsilon}\int_{x\in\Omega}\int_{z\in\Omega}\left(\Psi_{m(\sigma)}(z)-\Psi_{m(\sigma)}(x)\right)_+^2\theta^+(\mathrm{d}z)\theta^+(\mathrm{d}x)+\frac{\chi}{\varepsilon}
\tag{81}
$$

and

$$m(\omega)^2 = \frac{1}{2\varepsilon}\int_{x\in\Omega}\int_{z\in\Omega}\left(\Psi_{m(\omega)}(z)-\Psi_{m(\omega)}(x)\right)_+^2 \vartheta^+(\mathrm{d}z)\,\vartheta^+(\mathrm{d}x)+\frac{\chi}{\varepsilon}. \tag{82}$$

Afterward, we obtain

$$\begin{aligned}
&\left|m(\sigma)^2 - m(\omega)^2\right| \\
&= \left|m(\sigma)-m(\omega)\right|\left|m(\sigma)+m(\omega)\right| \\
&\leq \frac{1}{2\varepsilon}\left|\begin{aligned}&\int_{x\in\Omega}\int_{z\in\Omega}\left(\Psi_{m(\sigma)}(z,\sigma)-\Psi_{m(\sigma)}(x,\sigma)\right)_+^2\theta^+(\mathrm{d}z)\,\theta^+(\mathrm{d}x)\\ &-\int_{x\in\Omega}\int_{z\in\Omega}\left(\Psi_{m(\omega)}(z,\omega)-\Psi_{m(\omega)}(x,\omega)\right)_+^2\vartheta^+(\mathrm{d}z)\,\vartheta^+(\mathrm{d}x)\end{aligned}\right| \\
&\leq \frac{1}{2\varepsilon}\left|\begin{aligned}&\int_{x\in\Omega}\int_{z\in\Omega}\left(\Psi_{m(\sigma)}(z,\sigma)-\Psi_{m(\sigma)}(x,\sigma)\right)_+^2\theta^+(\mathrm{d}z)\,\theta^+(\mathrm{d}x)\\ &-\int_{x\in\Omega}\int_{z\in\Omega}\left(\Psi_{m(\sigma)}(z,\sigma)-\Psi_{m(\sigma)}(x,\sigma)\right)_+^2\vartheta^+(\mathrm{d}z)\,\vartheta^+(\mathrm{d}x)\end{aligned}\right| \\
&\quad+\frac{1}{2\varepsilon}\left|\begin{aligned}&\int_{x\in\Omega}\int_{z\in\Omega}\left(\Psi_{m(\sigma)}(z,\sigma)-\Psi_{m(\sigma)}(x,\sigma)\right)_+^2\vartheta^+(\mathrm{d}z)\,\vartheta^+(\mathrm{d}x)\\ &-\int_{x\in\Omega}\int_{z\in\Omega}\left(\Psi_{m(\omega)}(z,\omega)-\Psi_{m(\omega)}(x,\omega)\right)_+^2\vartheta^+(\mathrm{d}z)\,\vartheta^+(\mathrm{d}x)\end{aligned}\right| \\
&\leq \frac{1}{2\varepsilon}\bar{U}^2\left|\int_{x\in\Omega}\int_{z\in\Omega}\left(\theta^+(\mathrm{d}z)-\vartheta^+(\mathrm{d}z)\right)\left(\theta^+(\mathrm{d}x)-\vartheta^+(\mathrm{d}x)\right)\right| \\
&\quad+\frac{1}{2\varepsilon}\int_{x\in\Omega}\int_{z\in\Omega}\left|\left(\Psi_{m(\sigma)}(z,\sigma)-\Psi_{m(\sigma)}(x,\sigma)\right)_+^2-\left(\Psi_{m(\omega)}(z,\omega)-\Psi_{m(\omega)}(x,\omega)\right)_+^2\right|\vartheta^+(\mathrm{d}z)\,\vartheta^+(\mathrm{d}x) \\
&\leq \frac{1}{2\varepsilon}\bar{U}^2\left(\left\|\theta-\vartheta\right\|_{\mathrm{TV}}\right)^2+\frac{2\bar{U}}{\varepsilon}\left\|\Psi_{m(\sigma)}(\cdot,\sigma)-\Psi_{m(\omega)}(\cdot,\omega)\right\|_\infty\int_{x\in\Omega}\int_{z\in\Omega}\vartheta^+(\mathrm{d}z)\,\vartheta^+(\mathrm{d}x) \\
&\leq \frac{1}{2\varepsilon}\bar{U}^2\left(\left\|\theta-\vartheta\right\|_{\mathrm{TV}}\right)^2+\frac{2\bar{U}}{\varepsilon}\left(\left\|\vartheta\right\|_{\mathrm{TV}}\right)^2\left\|\Psi_{m(\sigma)}(\cdot,\sigma)-\Psi_{m(\omega)}(\cdot,\omega)\right\|_\infty \\
&\leq \frac{1}{2\varepsilon}\bar{U}^2\left(\left\|\sigma\right\|_{\mathrm{TV}}+\left\|\omega\right\|_{\mathrm{TV}}\right)\left\|\sigma-\omega\right\|_{\mathrm{TV}}+\frac{8\bar{U}}{\varepsilon}\left\|\Psi_{m(\sigma)}(\cdot,\sigma)-\Psi_{m(\omega)}(\cdot,\omega)\right\|_\infty \\
&\leq \frac{2\bar{U}^2}{\varepsilon}\left\|\sigma-\omega\right\|_{\mathrm{TV}}+\frac{8\bar{U}}{\varepsilon}\left\|\Psi_{m(\sigma)}(\cdot,\sigma)-\Psi_{m(\omega)}(\cdot,\omega)\right\|_\infty .
\end{aligned} \tag{83}$$

Consequently, we have

$$\left|m(\sigma)-m(\omega)\right| \leq \frac{1}{\left|m(\sigma)+m(\omega)\right|}\left(\frac{2\bar{U}^2}{\varepsilon}\left\|\sigma-\omega\right\|_{\mathrm{TV}}+\frac{8\bar{U}}{\varepsilon}\left\|\Psi_{m(\sigma)}(\cdot,\sigma)-\Psi_{m(\omega)}(\cdot,\omega)\right\|_\infty\right). \tag{84}$$

***Step 5. Auxiliary HJB equation: Lipschitz continuity 2***

By (80) (with $m = m(\sigma)$ and $n = m(\omega)$), (84) and the strict bound (75), if $\delta > 4\underline{\eta}^{-1}\bar{U}$, then we obtain

$$\left\|\Psi_{m(\sigma)}(\cdot,\sigma)-\Psi_{m(\omega)}(\cdot,\omega)\right\|_\infty \leq \left(1-\frac{4\bar{U}}{\delta\underline{\eta}}\right)^{-1}\left\{\left(L_U+\frac{\bar{U}^2}{2\delta\underline{\eta}}\right)\left\|\sigma-\omega\right\|_{\mathrm{TV}}+\frac{2\bar{U}^2}{\delta\underline{\eta}^2}\left|m(\sigma)-m(\omega)\right|\right\} \tag{85}$$

and

$$\left|m(\sigma)-m(\omega)\right| \leq \frac{1}{\underline{\eta}}\left(\frac{\bar{U}^2}{\varepsilon}\left\|\sigma-\omega\right\|_{\mathrm{TV}}+\frac{4\bar{U}}{\varepsilon}\left\|\Psi_{m(\sigma)}(\cdot,\sigma)-\Psi_{m(\omega)}(\cdot,\omega)\right\|_\infty\right). \tag{86}$$

Hence, we obtain

$$
\begin{aligned}
&\left\|\Psi_{m(\sigma)}(\cdot,\sigma)-\Psi_{m(\omega)}(\cdot,\omega)\right\|_{\infty} \\
&\leq\left(1-\frac{4\bar{U}}{\delta\underline{\eta}}\right)^{-1}\left(\begin{array}{l}\left(L_U+\dfrac{\bar{U}^2}{2\delta\underline{\eta}}\right)\|\sigma-\omega\|_{\mathrm{TV}} \\ +\dfrac{2\bar{U}^2}{\delta\underline{\eta}^3}\left(\dfrac{\bar{U}^2}{\varepsilon}\|\sigma-\omega\|_{\mathrm{TV}}+\dfrac{4\bar{U}}{\varepsilon}\left\|\Psi_{m(\sigma)}(\cdot,\sigma)-\Psi_{m(\omega)}(\cdot,\omega)\right\|_{\infty}\right)\end{array}\right) \\
&=\left(1-\frac{4\bar{U}}{\delta\underline{\eta}}\right)^{-1}\left(\left(L_U+\frac{\bar{U}^2}{2\delta\underline{\eta}}+\frac{4\bar{U}^4}{\delta\underline{\eta}^3\varepsilon}\right)\|\sigma-\omega\|_{\mathrm{TV}}+\frac{8\bar{U}^3}{\delta\underline{\eta}^3\varepsilon}\left\|\Psi_{m(\sigma)}(\cdot,\sigma)-\Psi_{m(\omega)}(\cdot,\omega)\right\|_{\infty}\right).
\end{aligned}
\tag{87}
$$

Therefore, if

$$
1-\left(1-\frac{4\bar{U}}{\delta\underline{\eta}}\right)^{-1}\frac{8\bar{U}^3}{\delta\underline{\eta}^3\varepsilon}>0 \quad \text{and hence} \quad \delta>\frac{4\bar{U}}{\underline{\eta}}+\frac{8\bar{U}^3}{\underline{\eta}^3\varepsilon}, \tag{88}
$$

then we obtain

$$
\left\|\Psi_{m(\sigma)}(\cdot,\sigma)-\Psi_{m(\omega)}(\cdot,\omega)\right\|_{\infty}\leq c_1\|\sigma-\omega\|_{\mathrm{TV}}, \tag{89}
$$

where $c_1>0$ is a constant depending only on $\delta,\bar{U},L_U,\underline{\eta},\varepsilon$. At the same time, by (86) and (89), we obtai

$$
\left|m(\sigma)-m(\omega)\right|\leq c_2\|\sigma-\omega\|_{\mathrm{TV}}. \tag{90}
$$

where $c_2>0$ is a constant depending only on $\delta,\bar{U},L_U,\underline{\eta},\varepsilon$.

***Step 6. EGD***

Now, we can identify $\eta(\sigma)=m(\sigma)$ and $\Phi(\cdot,\sigma)=\Psi_{\eta(\sigma)}(\cdot,\sigma)$ for any $\sigma\in\mathcal{M}_2$. By **Steps 4-5**, we obtain the following inequality to show the desired Lipschitz continuity of $\eta^{-1}\Phi$: Fix a sufficiently large $\delta>0$ such that (88); then, for any $\sigma,\omega\in\mathcal{M}_2$ and $y\in\Omega$,

$$
\begin{aligned}
\left|\frac{\Phi(y,\sigma)}{\eta(\sigma)}-\frac{\Phi(y,\omega)}{\eta(\omega)}\right| &\leq\left|\frac{\Phi(y,\sigma)}{\eta(\sigma)}-\frac{\Phi(y,\sigma)}{\eta(\omega)}\right|+\left|\frac{\Phi(y,\sigma)}{\eta(\omega)}-\frac{\Phi(y,\omega)}{\eta(\omega)}\right| \\
&\leq\frac{\Phi(y,\sigma)}{\underline{\eta}^2}\left|\eta(\sigma)-\eta(\omega)\right|+\frac{1}{\underline{\eta}}\left|\Phi(y,\sigma)-\Phi(y,\omega)\right| \\
&\leq\frac{\bar{U}}{\underline{\eta}^2}\left|\eta(\sigma)-\eta(\omega)\right|+\frac{1}{\underline{\eta}}\left|\Phi(y,\sigma)-\Phi(y,\omega)\right| \\
&\leq\left(\frac{1}{\underline{\eta}}c_1+\frac{\bar{U}}{\underline{\eta}^2}c_2\right)\|\sigma-\omega\|_{\mathrm{TV}}.
\end{aligned}
\tag{91}
$$

The boundedness of $\eta^{-1}\Phi$ is clear because $0\leq\eta^{-1}\Phi\leq\underline{\eta}^{-1}\bar{U}$. Then, Theorem 1 in Cheung [12] is applied, and the proof of the proposition is completed.

□

***Proof of Proposition 2***

First, in view of Theorem 1 in Cheung [12], the main task is to show that the quantity $\eta(\sigma)^{-1}\Phi(\cdot,\sigma)$ at each $x\in\Omega$ is bounded and Lipschitz continuous with respect to an arbitrary $\sigma\in\mathcal{M}_2$. Additionally, we need to show the unique existence of the couple $(\eta(\sigma),\Phi(\cdot,\sigma))\in\mathbb{R}_+\times C$ for each $\sigma\in\mathcal{M}_2$ and their uniform bounds. We also obtain the estimate $\underline{\eta}\leq\eta(\cdot)\leq\overline{\eta}$. The proof is divided into several parts.

***Step 1. Auxiliary HJB equation: existence***

We fix arbitrary $m>0$ and $\sigma\in\mathcal{M}_2$. We show that the following auxiliary HJB equation yields a unique solution $\Psi\in C$ such that $0\leq\Psi\leq\overline{U}$ when $\delta$ is large:

$$\delta\Psi(x)=\delta U(x,\sigma)+\sup_{u_x(\cdot)\geq 0,\ \int_\Omega u_x(z)\mathrm{d}z=1}\left\{\begin{array}{l}\int_\Omega u_x(z)\left(\Psi(z)-\Psi(x)\right)\mathrm{d}z\\ -m\int_\Omega\left\{u_x(z)\ln\left(u_x(z)\right)-u_x(z)+1\right\}\mathrm{d}z\end{array}\right\},\quad x\in\Omega. \tag{92}$$

First, if $\Psi\in C$, then the second term on the right-hand side of (92) is rewritten as follows:

$$\sup_{u_x(\cdot)\geq 0,\ \int_\Omega u_x(z)\mathrm{d}z=1}\left\{\begin{array}{l}\int_\Omega u_x(z)\left(\Psi(z)-\Psi(x)\right)\mathrm{d}z\\ -m\int_\Omega\left\{u_x(z)\ln\left(u_x(z)\right)-u_x(z)+1\right\}\mathrm{d}z\end{array}\right\}=m\ln\left(\int_\Omega\exp\left(\frac{\Psi(z)}{m}\right)\mathrm{d}z\right)-\Psi(x) \tag{93}$$

because an elementary calculation shows that (actually, $u_x$ does not depend on $x$)

$$\underset{u_x(\cdot)\geq 0,\ \int_\Omega u_x(z)\mathrm{d}z=1}{\arg\max}\left\{\begin{array}{l}\int_\Omega u_x(z)\left(\Psi(z)-\Psi(x)\right)\mathrm{d}z\\ -m\int_\Omega\left\{u_x(z)\ln\left(u_x(z)\right)-u_x(z)+1\right\}\mathrm{d}z\end{array}\right\}=\frac{\exp\left(\frac{\Psi(\cdot)}{m}\right)}{\int_\Omega\exp\left(\frac{\Psi(z)}{m}\right)\mathrm{d}z}\quad\left(:=\hat{u}(\cdot)\right). \tag{94}$$

Then, the auxiliary HJB equation (92) is rewritten as follows:

$$\delta\Psi(x)=\delta U(x,\sigma)+m\ln\left(\int_\Omega\exp\left(\frac{\Psi(z)}{m}\right)\mathrm{d}z\right)-\Psi(x),\quad x\in\Omega, \tag{95}$$

and hence

$$\Psi(x)=\frac{\delta}{\delta+1}U(x,\sigma)+\frac{m}{\delta+1}\ln\left(\int_\Omega\exp\left(\frac{\Psi(z)}{m}\right)\mathrm{d}z\right)\quad\left(:=\mathbb{B}(\Psi)(x)\right),\quad x\in\Omega. \tag{96}$$

From (96), we obtain

$$\begin{aligned}\exp\left(\frac{\Psi(x)}{m}\right)&=\exp\left(\frac{\delta}{\delta+1}\frac{U(x,\sigma)}{m}+\frac{1}{\delta+1}\ln\left(\int_\Omega\exp\left(\frac{\Psi(z)}{m}\right)\mathrm{d}z\right)\right)\\ &=\left(\int_\Omega\exp\left(\frac{\Psi(z)}{m}\right)\mathrm{d}z\right)^{\frac{1}{\delta+1}}\exp\left(\frac{\delta}{\delta+1}\frac{U(x,\sigma)}{m}\right)\end{aligned},\quad x\in\Omega \tag{97}$$

and hence

$$\int_\Omega\exp\left(\frac{\Psi(z)}{m}\right)\mathrm{d}z=\left(\int_\Omega\exp\left(\frac{\delta}{\delta+1}\frac{U(z,\sigma)}{m}\right)\mathrm{d}z\right)^{\frac{\delta+1}{\delta}}. \tag{98}$$

Substituting (98) into (96) yields the closed-form solution to the auxiliary HJB equation (92):

$$\Psi(x)=\Psi_m(x):=\frac{\delta}{\delta+1}U(x,\sigma)+\frac{m}{\delta}\ln\left(\int_\Omega \exp\left(\frac{\delta}{\delta+1}\frac{U(z,\sigma)}{m}\right)\mathrm{d}z\right),\quad x\in\Omega . \tag{99}$$

Here, the right-hand side of (99) does not depend on $\Psi$. This solution is bounded below by 0:

$$\Psi_m(x)\geq\frac{\delta}{\delta+1}\cdot 0+\frac{m}{\delta}\ln\left(\int_\Omega \exp\left(\frac{\delta}{\delta+1}\frac{0}{m}\right)\mathrm{d}z\right)=0,\quad x\in\Omega . \tag{100}$$

Similarly, this solution is bounded above by $\bar{U}$:

$$\begin{aligned}\Psi_m(x)&\leq\frac{\delta}{\delta+1}\bar{U}+\frac{m}{\delta}\ln\left(\int_\Omega \exp\left(\frac{\delta}{\delta+1}\frac{\bar{U}}{m}\right)\mathrm{d}z\right)\\&=\frac{\delta}{\delta+1}\bar{U}+\frac{m}{\delta}\ln\left(\exp\left(\frac{\delta}{\delta+1}\frac{\bar{U}}{m}\right)\right)\\&=\bar{U}\end{aligned},\quad x\in\Omega . \tag{101}$$

Finally, $\Psi_m$ is not a constant function because of **Assumption 1**.

***Step 2. Auxiliary HJB equation: uniqueness***

We found a closed-form solution (99) to the auxiliary HJB equation (96) and hence proved its existence. Now, we move to uniqueness. To show this, it suffices to prove that $\mathbb{B}$ is bounded and strictly contractive in $C$ for a sufficiently large $\delta>0$. Note that $\mathbb{B}(\Psi)\in C$ when $\Psi\in C$ because of the continuity of $U$ in **(A4) in Assumption 1**, i.e., (19).

Under the setting of **Step 1**, any solution $\Psi\in C$ to the auxiliary HJB equation (96) is bounded as $0\leq\Psi\leq\bar{U}$ in $\Omega$. Indeed, if $\Psi$ is minimized at some $\hat{x}\in\Omega$, then

$$\Psi(\hat{x})\geq\frac{\delta}{\delta+1}\cdot 0+\frac{m}{\delta}\ln\left(\int_\Omega \exp\left(\frac{\delta}{\delta+1}\frac{\Psi(\hat{x})}{m}\right)\mathrm{d}z\right)=\frac{\Psi(\hat{x})}{\delta+1}, \tag{102}$$

yielding $\Psi(\hat{x})\geq 0$. Similarly, if $\Psi$ is maximized at some $\hat{x}\in\Omega$, then

$$\Psi(\hat{x})\leq\frac{\delta}{\delta+1}\bar{U}+\frac{m}{\delta}\ln\left(\int_\Omega \exp\left(\frac{\delta}{\delta+1}\frac{\Psi(\hat{x})}{m}\right)\mathrm{d}z\right)=\frac{\delta}{\delta+1}\bar{U}+\frac{\Psi(\hat{x})}{\delta+1}, \tag{103}$$

yielding $\Psi(\hat{x})\leq\bar{U}$.

The *a priori* bound derived above implies that for any solutions $\Psi\in C$ to the auxiliary HJB equation (96), it suffices to consider those with $0\leq\Psi\leq\bar{U}$ in $\Omega$, as in **Proof of Proposition 1**. We prove the boundedness and strict contraction property of $\mathbb{B}:C\to C$ if $\delta>0$ is large. The boundedness is proven as follows: For any $\Psi_1\in C$ with $\|\Psi_1\|_\infty\leq\bar{U}$, as for (101),

$$\|\mathbb{B}(\Psi_1)\|_\infty\leq\bar{U} . \tag{104}$$

The strict contraction property is proven as follows: For any $\Psi_1,\Psi_2\in C$ with $\|\Psi_1\|_\infty,\|\Psi_2\|_\infty\leq\bar{U}$, if $\delta>\exp\left(m^{-1}\bar{U}\right)-1$, then

$$
\begin{aligned}
\left\|\mathbb{B}(\Psi_1)-\mathbb{B}(\Psi_2)\right\|_\infty &= \frac{m}{\delta+1}\left|\ln\left(\int_\Omega \exp\left(m^{-1}\Psi_1(z)\right)\mathrm{d}z\right)-\ln\left(\int_\Omega \exp\left(m^{-1}\Psi_2(z)\right)\mathrm{d}z\right)\right| \\
&\le \frac{m}{\delta+1}\left|\int_\Omega \exp\left(m^{-1}\Psi_1(z)\right)\mathrm{d}z-\int_\Omega \exp\left(m^{-1}\Psi_2(z)\right)\mathrm{d}z\right| \\
&\le \frac{m}{\delta+1}\int_\Omega \left|\exp\left(m^{-1}\Psi_1(z)\right)-\exp\left(m^{-1}\Psi_2(z)\right)\right|\mathrm{d}z \\
&\le \frac{1}{\delta+1}\exp\left(m^{-1}\bar{U}\right)\left\|\Psi_1-\Psi_2\right\|_\infty \\
&< \left\|\Psi_1-\Psi_2\right\|_\infty .
\end{aligned}
\tag{105}
$$

Therefore, by Banach's fixed-point theorem (e.g., Theorem 5.7 in Brezis [9]), if $\delta > \exp\left(m^{-1}\bar{U}\right)-1$, then the auxiliary HJB equation admits the unique solution $\Psi_m$ of (99).

***Step 3. Auxiliary HJB equation: cost constraint***

From (97)-(99), we have

$$
\frac{\exp\left(m^{-1}\Psi_m(z)\right)}{\int_\Omega \exp\left(m^{-1}\Psi_m(z)\right)\mathrm{d}z}=\frac{\exp\left(m^{-1}W(z,\sigma)\right)}{\int_\Omega \exp\left(m^{-1}W(z,\sigma)\right)\mathrm{d}z},
\tag{106}
$$

where $W=\delta U/(\delta+1)$. Then, we obtain

$$
\begin{aligned}
&\int_{x\in\Omega}\int_{z\in\Omega}\left\{\hat{u}(z)\ln\left(\hat{u}(z)\right)-\hat{u}(z)+1\right\}\mathrm{d}z\mathrm{d}x \\
&=\int_{z\in\Omega}\left\{\hat{u}(z)\ln\left(\hat{u}(z)\right)-\hat{u}(z)+1\right\}\mathrm{d}z \\
&=\int_{z\in\Omega}\left(\frac{\exp\left(m^{-1}\Psi_m(z)\right)}{\int_\Omega \exp\left(m^{-1}\Psi_m(y)\right)\mathrm{d}y}\ln\left(\frac{\exp\left(m^{-1}\Psi_m(z)\right)}{\int_\Omega \exp\left(m^{-1}\Psi_m(y)\right)\mathrm{d}y}\right)\right)\mathrm{d}z \\
&=\int_{z\in\Omega}\left(\frac{\exp\left(m^{-1}W(z,\sigma)\right)}{\int_\Omega \exp\left(m^{-1}W(y,\sigma)\right)\mathrm{d}y}\ln\frac{\exp\left(m^{-1}W(z,\sigma)\right)}{\int_\Omega \exp\left(m^{-1}W(y,\sigma)\right)\mathrm{d}y}\right)\mathrm{d}z \\
&\left(:=g(m)\right).
\end{aligned}
\tag{107}
$$

For each $\sigma\in\mathcal{M}_2$, we show that exactly one $m=m(\sigma)>0$ exists such that

$$
g(m)=\varepsilon .
\tag{108}
$$

If such an $m(\sigma)>0$ exists, then the inequality (33) is satisfied with an equality, and $m(\sigma)$ is the optimal Lagrangian multiplier due to the classical KKT result because both the objective and constraint are concave (whose negative is convex) and the constraint satisfies the Slater condition (Theorem 3.9 in Barbu and Precupanu [5]). Differentiating $g(m)$ with respect to $m>0$ yields (the existence of each integral follows from the dominate convergence because of the boundedness of $W$ )

$$
\begin{aligned}
\frac{\mathrm{d}g(m)}{\mathrm{d}m} &= -\frac{1}{m^3}\left(\frac{\int_\Omega \left(W(y,\sigma)\right)^2\exp\left(m^{-1}W(y,\sigma)\right)\mathrm{d}y}{\int_\Omega \exp\left(m^{-1}W(y,\sigma)\right)\mathrm{d}y}-\left(\frac{\int_\Omega W(y,\sigma)\exp\left(m^{-1}W(y,\sigma)\right)\mathrm{d}y}{\int_\Omega \exp\left(m^{-1}W(y,\sigma)\right)\mathrm{d}y}\right)^2\right) \\
&= -\frac{1}{m^3}\mathbb{V}_{m,\sigma}\left[W(y,\sigma)\right] \\
&< 0,
\end{aligned}
\tag{109}
$$

where the variance $\mathbb{V}_{m,\sigma}$ is based on the probability density $p_{m,\sigma}(y)=\dfrac{\exp\left(m^{-1}W(y,\sigma)\right)}{\int_{\Omega}\exp\left(m^{-1}W(z,\sigma)\right)\mathrm{d}z}$. The last strict inequality comes from **(A1) in Assumption 1**, where we assume a nonconstant $U$. Moreover, $\mathbb{V}_{m,\sigma}\left[W(y,\sigma)\right]$ is bounded for each $m>0$ because $W$ is nonnegative and strictly bounded.

From (109), it follows that $g$ is strictly decreasing for $m>0$. Moreover, $\lim_{m\to+0}g(m)=0$ because $p_{m,\sigma}$ reduces to the probability density of the uniform distribution on $\Omega$. Moreover, the distribution $p_{m,\sigma}$ approaches a (superposition of) Dirac delta concentrated at a nonzero maximum point of $\exp\left(m^{-1}W(y,\sigma)\right)$ (e.g., Corollary 2.1 in Hwang [22]); hence, $\lim_{m\to+0}g(m)=+\infty$. This combined with (109) and $\lim_{m\to+0}g(m)=0$ shows that exactly one $m(\sigma)>0$ such that $g\left(m(\sigma)\right)=\varepsilon$ exists.

Finally, we obtain the upper and lower bounds of $m(\sigma)$. By (109), for the upper bound, we obtain (we set $\overline{W}=\delta\overline{U}/(\delta+1)$)

$$
\begin{aligned}
\varepsilon &= g\left(m(\sigma)\right)\\
&= \int_{m(\sigma)}^{+\infty}\frac{1}{n^{3}}\mathbb{V}_{n,\sigma}\left[W(y,\sigma)\right]\mathrm{d}n\\
&\leq \int_{m(\sigma)}^{+\infty}\frac{1}{n^{3}}\frac{\int_{\Omega}\left(W(y,\sigma)\right)^{2}\exp\left(n^{-1}W(y,\sigma)\right)\mathrm{d}y}{\int_{\Omega}\exp\left(n^{-1}W(y,\sigma)\right)\mathrm{d}y}\mathrm{d}n\\
&\leq \overline{W}^{2}\int_{m(\sigma)}^{+\infty}\frac{1}{n^{3}}\mathrm{d}n\\
&\leq \frac{\overline{W}^{2}}{2\left(m(\sigma)\right)^{2}},
\end{aligned}
\tag{110}
$$

and hence:

$$
m(\sigma)\leq\frac{\overline{W}}{\sqrt{2\varepsilon}}\leq\frac{\overline{U}}{\sqrt{2\varepsilon}}\left(:=\overline{\eta}\right).
\tag{111}
$$

For the lower bound, by **(A5) in Assumption 1** and

$$
\varepsilon=\int_{z\in\Omega}\left(\frac{\exp\left(n(\sigma)^{-1}U(z,\sigma)\right)}{\int_{\Omega}\exp\left(n(\sigma)^{-1}U(y,\sigma)\right)\mathrm{d}y}\ln\frac{\exp\left(n(\sigma)^{-1}U(z,\sigma)\right)}{\int_{\Omega}\exp\left(n(\sigma)^{-1}U(y,\sigma)\right)\mathrm{d}y}\right)\mathrm{d}z
\tag{112}
$$

with $n(\sigma)=\dfrac{\delta+1}{\delta}m(\sigma)$, we obtain

$$
n(\sigma)\geq\underline{\eta}_{\varepsilon}\leftrightarrow m(\sigma)\geq\frac{\delta}{\delta+1}\underline{\eta}_{\varepsilon}\left(:=\underline{\eta}\right).
\tag{113}
$$

Consequently, we obtain the bounds of $m(\sigma)$ independent of $\sigma\in\mathcal{M}_{2}$:

$$
\underline{\eta}\leq m(\sigma)\leq\overline{\eta}.
\tag{114}
$$

***Step 4. Auxiliary HJB equation: Lipschitz continuity 1***

For each $m>0$ and $\sigma \in \mathcal{M}_2$, we write the solution $\Psi_m(\cdot)$ in (99) as $\Psi_m(\cdot,\sigma)$. We show its Lipschitz continuity with respect to the second argument. Fixing $m,n>0$ and $\sigma,\omega \in \mathcal{M}_2$, we have

$$
\begin{aligned}
\left\|\Psi_m(\cdot,\sigma)-\Psi_n(\cdot,\omega)\right\|_\infty &\leq \frac{\delta}{\delta+1}\left\|U(\cdot,\sigma)-U(\cdot,\omega)\right\|_\infty \\
&\quad +\frac{1}{\delta}\left| m\ln\left(\int_\Omega \exp\left(\frac{\delta}{\delta+1}\frac{U(z,\sigma)}{m}\right)\mathrm{d}z\right)-n\ln\left(\int_\Omega \exp\left(\frac{\delta}{\delta+1}\frac{U(z,\omega)}{m}\right)\mathrm{d}z\right)\right| \\
&\leq \frac{\delta}{\delta+1}L_U\left\|\sigma-\omega\right\|_{\mathrm{TV}} \\
&\quad +\frac{1}{\delta}\left| m\ln\left(\int_\Omega \exp\left(m^{-1}W(z,\omega)\right)\mathrm{d}z\right)-n\ln\left(\int_\Omega \exp\left(m^{-1}W(z,\omega)\right)\mathrm{d}z\right)\right| \\
&\quad +\frac{1}{\delta}\left| m\ln\left(\int_\Omega \exp\left(m^{-1}W(z,\sigma)\right)\mathrm{d}z\right)-m\ln\left(\int_\Omega \exp\left(m^{-1}W(z,\omega)\right)\mathrm{d}z\right)\right| \\
&\leq \frac{\delta}{\delta+1}L_U\left\|\sigma-\omega\right\|_{\mathrm{TV}}+\frac{|m-n|}{\delta}\left|\ln\left(\int_\Omega \exp\left(m^{-1}W(z,\omega)\right)\mathrm{d}z\right)\right| \\
&\quad +\frac{m}{\delta}\left|\ln\left(\int_\Omega \exp\left(m^{-1}W(z,\sigma)\right)\mathrm{d}z\right)-\ln\left(\int_\Omega \exp\left(m^{-1}W(z,\omega)\right)\mathrm{d}z\right)\right| \\
&= \frac{\delta}{\delta+1}L_U\left\|\sigma-\omega\right\|_{\mathrm{TV}}+\frac{\bar{U}}{\delta+1}\frac{|m-n|}{m} \\
&\quad +\frac{m}{\delta}\left|\ln\left(\int_\Omega \exp\left(m^{-1}W(z,\sigma)\right)\mathrm{d}z\right)-\ln\left(\int_\Omega \exp\left(m^{-1}W(z,\omega)\right)\mathrm{d}z\right)\right|.
\end{aligned}
\tag{115}
$$

The last line of (115) is evaluated as follows:

$$
\begin{aligned}
&\left|\ln\left(\int_\Omega \exp\left(m^{-1}W(z,\sigma)\right)\mathrm{d}z\right)-\ln\left(\int_\Omega \exp\left(m^{-1}W(z,\omega)\right)\mathrm{d}z\right)\right| \\
&\leq \left|\int_\Omega \exp\left(m^{-1}W(z,\sigma)\right)\mathrm{d}z-\int_\Omega \exp\left(m^{-1}W(z,\omega)\right)\mathrm{d}z\right| \\
&\leq \left\|\exp\left(m^{-1}W(\cdot,\sigma)\right)-\exp\left(m^{-1}W(\cdot,\omega)\right)\right\|_\infty \\
&\leq \frac{\exp\left(m^{-1}\bar{W}\right)}{m}\left\|W(\cdot,\sigma)-W(\cdot,\omega)\right\|_\infty \\
&\leq \frac{\delta^2 L_U}{(\delta+1)^2 m}\exp\left(\frac{\delta}{\delta+1}\frac{\bar{U}}{m}\right)\left\|\sigma-\omega\right\|_{\mathrm{TV}}.
\end{aligned}
\tag{116}
$$

Substituting (116) into (115) yields

$$
\begin{aligned}
&\left\|\Psi_m(\cdot,\sigma)-\Psi_n(\cdot,\omega)\right\|_\infty \\
&\leq \frac{\delta}{\delta+1}L_U\left\|\sigma-\omega\right\|_{\mathrm{TV}}+\frac{\bar{U}}{\delta+1}\frac{|m-n|}{m}+\frac{m}{\delta}\frac{\delta^2 L_U}{(\delta+1)^2 m}\exp\left(\frac{\delta}{\delta+1}\frac{\bar{U}}{m}\right)\left\|\sigma-\omega\right\|_{\mathrm{TV}} \\
&= \frac{\delta L_U}{\delta+1}\left\{1+\frac{1}{\delta+1}\exp\left(\frac{\delta}{\delta+1}\frac{\bar{U}}{m}\right)\right\}\left\|\sigma-\omega\right\|_{\mathrm{TV}}+\frac{\bar{U}}{\delta+1}\frac{|m-n|}{m}.
\end{aligned}
\tag{117}
$$

***Step 5. Optimal Lagrangian multiplier: Lipschitz continuity***

According to **Step 3**, for each $\sigma \in \mathcal{M}_2$, there exists a unique $m=m(\sigma)>0$ such that $g(m(\sigma))=\varepsilon$. We evaluate the difference between $m(\sigma)$ and $m(\omega)$ for any $\sigma,\omega \in \mathcal{M}_2$. Below, we explicitly write the dependence of $g$ on $\sigma$ (resp., $\omega$) as $g_\sigma$ (resp., $g_\omega$).

By the classical mean value theorem, we have

$$\left|g_\sigma\left(m(\sigma)\right)-g_\sigma\left(m(\omega)\right)\right|=\left|\frac{\mathrm{d}g_\sigma\left(\hat{m}\right)}{\mathrm{d}m}\right|\left|m(\sigma)-m(\omega)\right|, \tag{118}$$

where $\hat{m}$ is between $m(\sigma)$ and $m(\omega)$. By (118), we obtain

$$\begin{aligned}\left|m(\sigma)-m(\omega)\right| &= \left|\frac{\mathrm{d}g_\sigma\left(\hat{m}\right)}{\mathrm{d}m}\right|^{-1}\left|g_\sigma\left(m(\sigma)\right)-g_\omega\left(m(\omega)\right)+g_\omega\left(m(\omega)\right)-g_\sigma\left(m(\omega)\right)\right| \\ &= \left|\frac{\mathrm{d}g_\sigma\left(\hat{m}\right)}{\mathrm{d}m}\right|^{-1}\left|g_\omega\left(m(\omega)\right)-g_\sigma\left(m(\omega)\right)\right| \\ &= \left|\frac{\mathrm{d}g_\sigma\left(\hat{m}\right)}{\mathrm{d}m}\right|^{-1}\left|\int_{z\in\Omega}p_{m,\omega}(z)\ln p_{m,\omega}(z)\mathrm{d}z-\int_{z\in\Omega}p_{m,\sigma}(z)\ln p_{m,\sigma}(z)\mathrm{d}z\right|_{m=m(\omega)}.\end{aligned} \tag{119}$$

The last absolute value is evaluated as follows:

$$\begin{aligned}&\left|\int_{z\in\Omega}p_{m,\omega}(z)\ln p_{m,\omega}(z)\mathrm{d}z-\int_{z\in\Omega}p_{m,\sigma}(z)\ln p_{m,\sigma}(z)\mathrm{d}z\right| \\ &\leq\left|\int_{z\in\Omega}p_{m,\omega}(z)\ln p_{m,\omega}(z)\mathrm{d}z-\int_{z\in\Omega}p_{m,\sigma}(z)\ln p_{m,\omega}(z)\mathrm{d}z\right| \\ &+\left|\int_{z\in\Omega}p_{m,\sigma}(z)\ln p_{m,\omega}(z)\mathrm{d}z-\int_{z\in\Omega}p_{m,\sigma}(z)\ln p_{m,\sigma}(z)\mathrm{d}z\right| \\ &\leq\int_{z\in\Omega}\left|p_{m,\omega}(z)-p_{m,\sigma}(z)\right|\left|\ln p_{m,\omega}(z)\right|\mathrm{d}z+\int_{z\in\Omega}p_{m,\sigma}(z)\left|\ln p_{m,\omega}(z)-\ln p_{m,\sigma}(z)\right|\mathrm{d}z.\end{aligned} \tag{120}$$

The first integral in the last line of (120) is evaluated as follows:

$$\begin{aligned}&\int_{z\in\Omega}\left|p_{m,\omega}(z)-p_{m,\sigma}(z)\right|\left|\ln p_{m,\omega}(z)\right|\mathrm{d}z \\ &=\int_{z\in\Omega}\left|p_{m,\omega}(z)-p_{m,\sigma}(z)\right|\left|\ln\frac{\exp\left(m^{-1}W(y,\omega)\right)}{\int_\Omega\exp\left(m^{-1}W(z,\omega)\right)\mathrm{d}z}\right|\mathrm{d}y \\ &\leq\int_{z\in\Omega}\left|p_{m,\omega}(z)-p_{m,\sigma}(z)\right|\left|\ln\exp\left(m^{-1}\bar{W}\right)\right|\mathrm{d}y \\ &\leq\frac{\bar{W}}{m}\int_{z\in\Omega}\left|p_{m,\omega}(z)-p_{m,\sigma}(z)\right|\mathrm{d}y \\ &\leq\frac{\bar{W}}{m}\frac{\int_{z\in\Omega}\left|\begin{array}{l}\exp\left(m^{-1}W(y,\omega)\right)\int_\Omega\exp\left(m^{-1}W(z,\sigma)\right)\mathrm{d}z \\ -\int_\Omega\exp\left(m^{-1}W(z,\omega)\right)\mathrm{d}z\exp\left(m^{-1}W(y,\sigma)\right)\end{array}\right|\mathrm{d}y}{\int_\Omega\exp\left(m^{-1}W(z,\omega)\right)\mathrm{d}z\int_\Omega\exp\left(m^{-1}W(z,\sigma)\right)\mathrm{d}z} \\ &\leq\frac{\bar{W}}{m}\int_{z\in\Omega}\left(\begin{array}{l}\int_\Omega\exp\left(m^{-1}W(z,\sigma)\right)\mathrm{d}z\left|\exp\left(m^{-1}W(y,\omega)\right)-\exp\left(m^{-1}W(y,\sigma)\right)\right| \\ +\exp\left(m^{-1}W(y,\sigma)\right)\left|\int_\Omega\exp\left(m^{-1}W(z,\sigma)\right)\mathrm{d}z-\int_\Omega\exp\left(m^{-1}W(z,\omega)\right)\mathrm{d}z\right|\end{array}\right)\mathrm{d}y \\ &\leq\frac{2\bar{W}}{m}\exp\left(m^{-1}\bar{W}\right)\int_{z\in\Omega}\left|\exp\left(m^{-1}W(z,\omega)\right)-\exp\left(m^{-1}W(z,\sigma)\right)\right|\mathrm{d}y \\ &\leq\frac{2\bar{W}}{m^2}\exp\left(2m^{-1}\bar{W}\right)\left\|W(\cdot,\omega)-W(\cdot,\sigma)\right\|_\infty \\ &\leq\frac{2\bar{W}L_U}{m^2}\exp\left(2m^{-1}\bar{W}\right)\left\|\omega-\sigma\right\|_{\mathrm{TV}}.\end{aligned} \tag{121}$$

Similarly, the second integral in the last line of (120) is evaluated as follows:

$$
\begin{aligned}
&\int_{z\in\Omega} p_{m,\sigma}(z)\left|\ln p_{m,\omega}(z)-\ln p_{m,\sigma}(z)\right|\mathrm{d}z\\
&=\int_{z\in\Omega}\frac{\exp\left(m^{-1}W(y,\sigma)\right)}{\int_{\Omega}\exp\left(m^{-1}W(z,\sigma)\right)\mathrm{d}z}\left|\ln p_{m,\omega}(z)-\ln p_{m,\sigma}(z)\right|\mathrm{d}y\\
&\le \exp\left(m^{-1}\bar{W}\right)\int_{z\in\Omega}\min\left\{p_{m,\omega}(z),p_{m,\sigma}(z)\right\}^{-1}\left|p_{m,\omega}(z)-p_{m,\sigma}(z)\right|\mathrm{d}y\\
&\le \exp\left(2m^{-1}\bar{W}\right)\int_{z\in\Omega}\left|p_{m,\omega}(z)-p_{m,\sigma}(z)\right|\mathrm{d}y\\
&\le \exp\left(2m^{-1}\bar{W}\right)\times\frac{m}{\bar{W}}\frac{2\bar{W}L_U}{m^2}\exp\left(2m^{-1}\bar{W}\right)\left\|\omega-\sigma\right\|_{\mathrm{TV}}\\
&=\frac{2L_U}{m}\exp\left(4m^{-1}\bar{W}\right)\left\|\omega-\sigma\right\|_{\mathrm{TV}}.
\end{aligned}
\tag{122}
$$

Consequently, we obtain

$$
\begin{aligned}
&\left|m(\sigma)-m(\omega)\right|\\
&\le\left|\frac{dg_\sigma(\hat{m})}{\mathrm{d}m}\right|^{-1}\left(\frac{2\bar{W}L_U}{m^2}\exp\left(2m^{-1}\bar{W}\right)\left\|\omega-\sigma\right\|_{\mathrm{TV}}+\frac{2L_U}{m}\exp\left(4m^{-1}\bar{W}\right)\left\|\omega-\sigma\right\|_{\mathrm{TV}}\right)_{m=m(\omega)}\\
&=\left|\frac{dg_\sigma(\hat{m})}{\mathrm{d}m}\right|^{-1}\left(\frac{\bar{W}}{m(\omega)}+\exp\left(2m(\omega)^{-1}\bar{W}\right)\right)\frac{2L_U}{m(\omega)}\exp\left(2m(\omega)^{-1}\bar{W}\right)\left\|\omega-\sigma\right\|_{\mathrm{TV}}.
\end{aligned}
\tag{123}
$$

***Step 6. Auxiliary HJB equation: Lipschitz continuity 2***

Fix s sufficiently large $\delta>0$ such that the following condition is satisfied (this is possible because the right-hand side of the inequality goes to $\exp\left(\underline{\eta}_{\varepsilon}^{-1}\bar{U}\right)-1$ as $\delta\to+\infty$):

$$
\delta>\exp\left(\underline{\eta}^{-1}\bar{U}\right)-1=\exp\left(\underline{\eta}_{\varepsilon}^{-1}\frac{(\delta+1)\bar{U}}{\delta}\right)-1.
\tag{124}
$$

By (117) and (123), substituting $m=m(\sigma)$ and $n=m(\omega)$ yields ($\hat{m}$ is between $m(\sigma)$ and $m(\omega)$, and hence $\underline{\eta}\le\hat{m}\le\bar{\eta}$ )

$$
\begin{aligned}
\left\|\Psi_{m(\sigma)}(\cdot,\sigma)-\Psi_{m(\omega)}(\cdot,\omega)\right\|_{\infty}&\le\frac{\delta L_U}{\delta+1}\left\{1+\frac{1}{\delta+1}\exp\left(\frac{\delta}{\delta+1}\frac{\bar{U}}{\underline{\eta}}\right)\right\}\left\|\sigma-\omega\right\|_{\mathrm{TV}}\\
&\quad+\frac{\bar{U}}{\delta+1}\frac{1}{\underline{\eta}}\left|\frac{dg_\sigma(\hat{m})}{\mathrm{d}m}\right|^{-1}\left(\frac{\bar{W}}{\underline{\eta}}+\exp\left(\frac{2\bar{W}}{\underline{\eta}}\right)\right)\frac{2L_U}{\underline{\eta}}\exp\left(\frac{2\bar{W}}{\underline{\eta}}\right)\left\|\omega-\sigma\right\|_{\mathrm{TV}}\\
&\le c_4\left\|\omega-\sigma\right\|_{\mathrm{TV}},
\end{aligned}
\tag{125}
$$

where $c_4>0$ is a constant depending only on $\delta,\bar{U},L_U,\underline{\eta},\bar{\eta},\varepsilon$. Here, we used the estimate

$$
\begin{aligned}
\left|\frac{dg_\sigma(\hat{m})}{\mathrm{d}m}\right|^{-1}&=\frac{\hat{m}^3}{\mathbb{V}_{\hat{m},\sigma}\left[W(y,\sigma)\right]}\\
&\le\frac{\bar{\eta}^3}{\mathbb{V}_{\hat{m},\sigma}\left[\hat{m}^{-1}W(y,\sigma)\right]}\\
&\le\frac{\bar{\eta}^3}{\underline{c}}
\end{aligned}
\tag{126}
$$

with a constant $\underline{c} > 0$ independent of $\sigma$. The reason that such a constant $\underline{c} > 0$ exists is as follows. First, we have

$$\frac{\int_\Omega (f(y))^2 \exp(f(y))\mathrm{d}y}{\int_\Omega \exp(f(y))\mathrm{d}y} - \left(\frac{\int_\Omega f(y)\exp(f(y))\mathrm{d}y}{\int_\Omega \exp(f(y))\mathrm{d}y}\right) > 0 \tag{127}$$

for any $f:\Omega \to \mathbb{R}_+$ that is Lipschitz continuous with some Lipschitz constant $A > 0$ and $\max_{x\in\Omega} f(x) - \min_{x\in\Omega} f(x) > B$ with some constant $B > 0$. The left-hand side of (127) is the variance of $f(x)$ with the probability density proportional to $\exp(f(x))$. The support of this density is $\Omega$, and $f$ is not a constant function. This means that the variance of $f(x)$ never becomes 0, and hence, the left-hand side of (127) is bounded from below by some constant $c(A,B) > 0$ depending only on $A, B$. Applying this observation along with **(A1) and (A4)** shows the existence of $\underline{c}$.

***Step 7. EGD***

Now, we can identify $\eta(\sigma) = m(\sigma)$ and $\Phi(\cdot,\sigma) = \Psi_{\eta(\sigma)}(\cdot,\sigma)$ for any $\sigma \in \mathcal{M}_2$. By **Steps 5-6**, we obtain the following inequality to show the desired Lipschitz continuity of $\eta^{-1}\Phi$: For any $\sigma, \omega \in \mathcal{M}_2$ and $y \in \Omega$,

$$\begin{aligned}
\left|\frac{\Phi(y,\sigma)}{\eta(\sigma)} - \frac{\Phi(y,\omega)}{\eta(\omega)}\right| &\leq \left|\frac{\Phi(y,\sigma)}{\eta(\sigma)} - \frac{\Phi(y,\sigma)}{\eta(\omega)}\right| + \left|\frac{\Phi(y,\sigma)}{\eta(\omega)} - \frac{\Phi(y,\omega)}{\eta(\omega)}\right| \\
&\leq \frac{\Phi(y,\sigma)}{\eta(\sigma)\eta(\omega)}\left|\eta(\sigma) - \eta(\omega)\right| + \frac{1}{\eta(\omega)}\left|\Phi(y,\sigma) - \Phi(y,\omega)\right| \\
&\leq \frac{\bar{U}}{\underline{\eta}^2}\left|\eta(\sigma) - \eta(\omega)\right| + \frac{1}{\underline{\eta}}\left|\Phi(y,\sigma) - \Phi(y,\omega)\right| \\
&\leq c_5 \left\|\sigma - \omega\right\|_{\mathrm{TV}}
\end{aligned} \tag{128}$$

with a constant $c_5 > 0$ independent of $\sigma, \omega \in \mathcal{M}_2$ and $y \in \Omega$. The boundedness of $\eta^{-1}\Phi \leq \underline{\eta}^{-1}\bar{U}$ is clear. Then, Theorem 1 in Cheung [12] is applied, and the proof of the proposition is completed.

□

### A.3 Two-dimensional numerical method

We explain how to discretize the logit model in the two-dimensional $x$-$z$ space. The explanation below is a direct extension of that for the numerical method presented in the main test. We use a finite difference method based on a fixed time increment $\Delta t > 0$ and spatial increment $\Delta x = \Delta z = 1/N$, where $N \in \mathbb{N}$ is the degree of freedom in the $x$ direction. We prepare one-dimensional cells $\omega_i = \left[(i-1)1/N, i/N\right)$ for $i = 1,2,...,N-1$ and $\omega_N = \left[1-1/N,1\right]$, and two-dimensional cells as $\Omega_{i,j} = \omega_i \times \omega_j$ ($i,j = 1,2,...,N$). Each quantity discretized at time step $k\Delta t$ ($k = 0,1,2,...$) and in $\Omega_{i,j}$ ($i,j = 1,2,...,N$) is denoted by the subscript $i,j,k$. The initial condition $\mu_0$ is discretized in each $\Omega_{i,j}$ as $\mu_{0,i,j} = \int_{\Omega_{i,j}} \mu_0(\mathrm{d}x)$ or its approximation if necessary.

Assume that we already have $\mu_{i,j,k-1}$ ($i,j = 1,2,...,N$) for some $k \in \mathbb{N}$. Each $\mu_{i,j,k}$ ($i,j = 1,2,...,N$) is computed by applying a fully explicit Euler discretization:

$$\mu_{i,j,k} = (1-\Delta t)\mu_{i,j,k-1} + \Delta t \frac{\exp\left(\dfrac{\Phi_{i,j,k-1}}{\eta_{k-1}}\right)}{\sum_{l,m=1}^{N} \exp\left(\dfrac{\Phi_{l,m,k-1}}{\eta_{k-1}}\right)\Delta x \Delta xy}. \tag{129}$$

The quantities $\Phi_{i,j,k-1}$ ($i,j = 1,2,...,N$) and $\eta_{k-1}$ are found in the following system, where $\Phi_{i,j,k-1}$ is found explicitly because of (35):

$$\Phi_{i,j,k-1} = U_{i,j,k-1} + \frac{\eta_{k-1}}{\delta} \ln\left(\sum_{l,m=1}^{N} \exp\left(\frac{\delta}{\delta+1}\frac{U_{l,m,k-1}}{\eta_{k-1}}\right)\Delta x \Delta z\right), \quad i,j = 1,2,...,N \tag{130}$$

and

$$\varepsilon = \sum_{i,j=1}^{N} \left\{ \frac{\exp\left(\dfrac{\delta}{\delta+1}\dfrac{U_{i,j,k-1}}{\eta_{k-1}}\right)}{\sum_{l,m=1}^{N} \exp\left(\dfrac{\delta}{\delta+1}\dfrac{U_{l,m,k-1}}{\eta_{k-1}}\right)\Delta x \Delta z} \ln\left( \frac{\exp\left(\dfrac{\delta}{\delta+1}\dfrac{U_{i,j,k-1}}{\eta_{k-1}}\right)}{\sum_{l,m=1}^{N} \exp\left(\dfrac{\delta}{\delta+1}\dfrac{U_{l,m,k-1}}{\eta_{k-1}}\right)\Delta x \Delta z} \right) \right\} \Delta x \Delta z \ \left(:= g_{\Delta}(\eta_{k-1})\right). \tag{131}$$

The system (131) is iteratively solved for $\eta_{k-1}$ by using **Algorithm 2** as in the one-dimensional case, where $g_{k-1}^{(n)}$ is set as

$$g_{k-1}^{(n)} = \frac{\delta}{\delta+1} \frac{1}{\varepsilon + \ln\left(\sum_{i,j=1}^{N} \exp\left(\dfrac{\delta}{\delta+1}\dfrac{U_{i,j,k-1}}{\eta_{k-1}^{(n)}}\right)\Delta x\right)} \frac{\sum_{i,j=1}^{N} U_{i,j,k-1} \exp\left(\dfrac{\delta}{\delta+1}\dfrac{U_{i,j,k-1}}{\eta_{k-1}^{(n)}}\right)\Delta x}{\sum_{i,j=1}^{N} \exp\left(\dfrac{\delta}{\delta+1}\dfrac{U_{i,j,k-1}}{\eta_{k-1}^{(n)}}\right)\Delta x}. \tag{132}$$

### A.4 Convergence study

This section briefly conducts a convergence study of the finite difference method to justify the computational resolution employed in **Section 4**. We present the result for the replicator dynamics here, but comparable results apply to the other two models. **Figure A1** compares stationary probability densities $p$ with $\Delta x = 1/250$ and $\Delta x = 1/500$ while $\Delta t = 0.005$ is fixed, demonstrating that they are visually indifferent; their absolute difference is estimated to be 0.0024% on average in $\Omega$. **Figure A2** compares the time histories of the optimal Lagrangian $\eta$ stationary probability densities with $\Delta t = 0.005$ and $\Delta x = 0.0025$ while $\Delta x = 1/250$ is fixed, again demonstrating that they are visually indifferent; their absolute difference is estimated to be on average 0.25% during $0 \le t \le 10$. Considering the conventional fact that a common finite difference method exhibits first- to second-order accuracy in space and time (e.g., Yoshioka and Tsujimura [50]), these computational results imply that the computational resolution $\Delta x = 1/250$ and $\Delta t = 0.005$ used in **Section 4** is fine enough for the discussion made in that section.

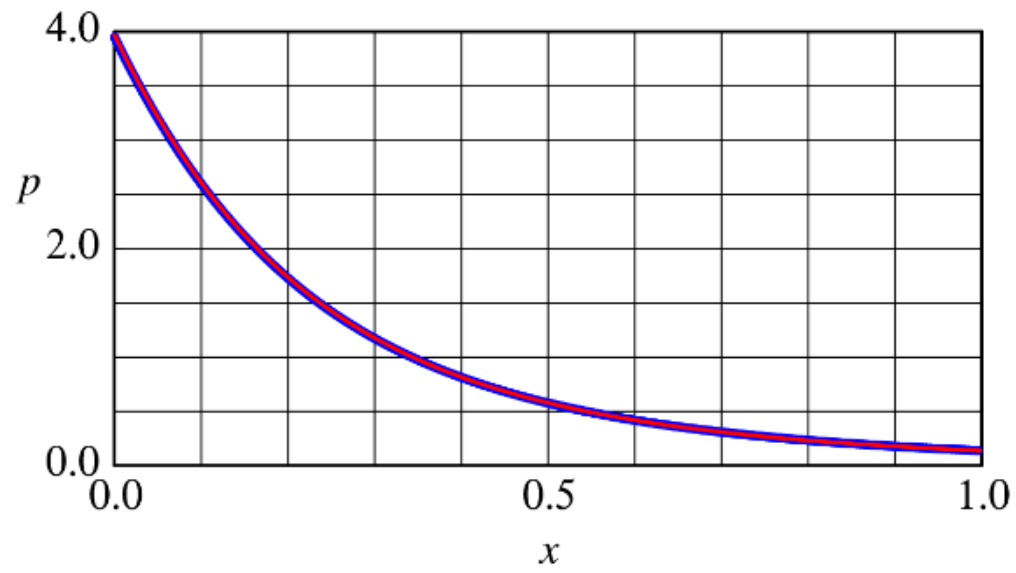

**Figure A1.** Comparison of the computed stationary probability densities $p$: $\Delta x = 1/250$ (blue) and $\Delta x = 1/500$ (red).

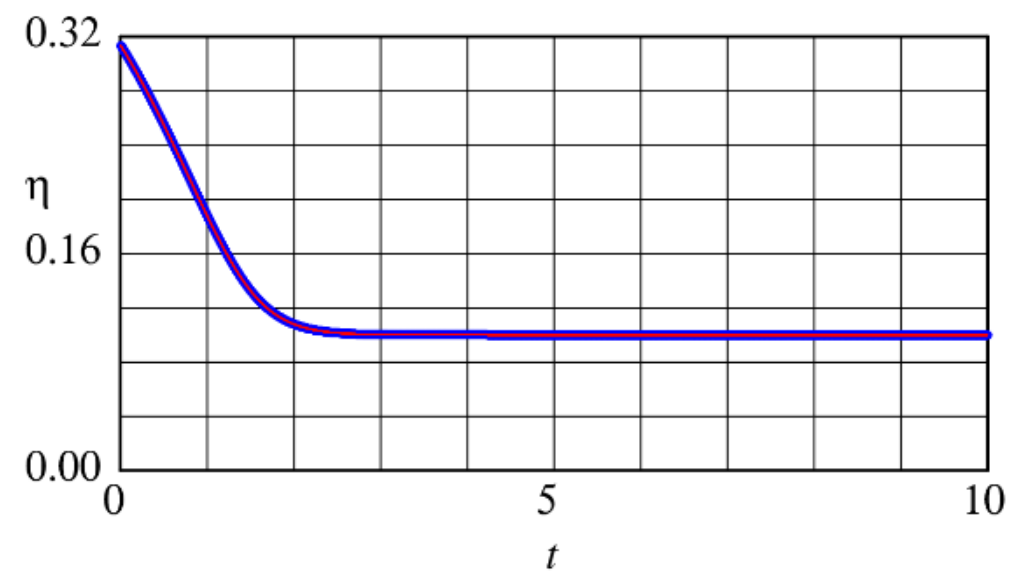

**Figure A2.** Comparison time histories of the computed optimal Lagrangian $\eta$: $\Delta t = 0.005$ (blue) and $\Delta t = 0.0025$ (red).

### A.5 Analytical results about a utility

For the following utility ( $\mathbb{E}_\mu$ and $\mathbb{V}_\mu$ represent expectation and variance with respect to $\mu \in \mathcal{P}$, respectively)

$$U_1(x,\mu) = \int_\Omega (x-y)^2 \mu(\mathrm{d}y) = \left(x - \mathbb{E}_\mu[x]\right)^2 + \mathbb{V}_\mu[x], \tag{133}$$

the cost constraint (33), which is satisfied at an equality because of **Proof of Proposition 2**, is reduced to

$$\varepsilon = \int_{z\in\Omega} \left( \frac{\exp\left(m^{-1}\left(y - \mathbb{E}_\mu[x]\right)^2\right)}{\int_\Omega \exp\left(m^{-1}\left(y - \mathbb{E}_\mu[x]\right)^2\right)\mathrm{d}y} \ln \frac{\exp\left(m^{-1}\left(y - \mathbb{E}_\mu[x]\right)^2\right)}{\int_\Omega \exp\left(m^{-1}\left(y - \mathbb{E}_\mu[x]\right)^2\right)\mathrm{d}y} \right) \mathrm{d}z\,. \tag{134}$$

If $\mu$ is symmetric, i.e., if $\mathbb{E}_\mu[x] = 1/2$, then we obtain

$$\varepsilon = \int_{z\in\Omega} \left( \frac{\exp\left(m^{-1}\left(y-1/2\right)^2\right)}{\int_{\Omega} \exp\left(m^{-1}\left(y-1/2\right)^2\right)\mathrm{d}y} \ln \frac{\exp\left(m^{-1}\left(y-1/2\right)^2\right)}{\int_{\Omega} \exp\left(m^{-1}\left(y-1/2\right)^2\right)\mathrm{d}y} \right) \mathrm{d}z \ . \tag{135}$$

Because the right-hand side of (135) is independent of $\mu$, this equation is simply an algebraic equation of $m > 0$ for each given $\varepsilon > 0$. This is the reason that the optimal Lagrangian multiplier $\eta_t$ is time-independent in **Figures 1 and 3** in **Section 4**. Our numerical solutions are therefore consistent with the theory.